\documentclass[11pt,a4paper]{article}
\usepackage[margin=1in]{geometry}
\usepackage{amsmath,amssymb,amsfonts,bm,mathtools}
\usepackage{graphicx}
\usepackage{subcaption}
\usepackage{booktabs}
\usepackage{multirow}
\usepackage{array}
\usepackage{enumitem}
\usepackage[numbers,sort&compress]{natbib}
\usepackage[hidelinks]{hyperref}
\usepackage{url}
\usepackage{amsthm}
\theoremstyle{plain}
\newtheorem{theorem}{Theorem}[section]
\newtheorem{proposition}[theorem]{Proposition}

\theoremstyle{definition}

\theoremstyle{remark}

\numberwithin{equation}{section}
\title{Spectral-Domain Information Geometry of Stationary Covariance Models:
	Gaussian and Elliptical Extensions}

\author{
	Shu Yang, Dongying Li, Wenxian Yu
}

\date{}

\begin{document}
	
	\maketitle
	
\begin{abstract}
This paper develops an information-geometric framework for stationary covariance structures. We begin with the infinite-dimensional setting, where stationary covariance structures can be characterized without finite-window boundary effects. By formulating the covariance structure as a bilinear form on a suitable class of test functions, we establish its spectral representation and show how finite-dimensional Toeplitz covariance matrices arise as restrictions to nested observation subspaces. For zero-mean proper complex Gaussian wide-sense stationary (WSS) processes, we show that the normalized finite-dimensional Fisher information converges to a spectral-domain metric. In logarithmic spectral coordinates, this metric becomes the standard $L^2$ inner product, yielding an explicit Fisher--Rao distance that is also recovered as the asymptotic normalized affine-invariant Riemannian distance between Toeplitz covariance matrices.

We then extend the framework to elliptical observation models by separating global covariance scale from shape. The shape Fisher information is orthogonal to radial nuisance directions and is characterized by a dimension- and radial-model-dependent coefficient. We derive a condition under which this coefficient converges to its Gaussian value as the observation dimension tends to infinity. The resulting limiting spectral-shape geometry is the centered $L^2$ geometry of log-spectra and is invariant to positive rescaling. Numerical experiments illustrate the theoretical convergence results and the recovery of the Gaussian shape coefficient.
	
	\end{abstract}
	

	\section{Introduction}
	
	Stationary random processes are commonly characterized through their second-order statistics. For a wide-sense stationary (WSS) process, the autocorrelation function and its spectral representation describe the dependence structure across time and frequency \cite{brockwell2009time,hannan2009multiple,ibragimov2012gaussian}. In practical signal processing, however, only finite segments of an underlying process are available for observation. The corresponding covariance matrices are finite-dimensional restrictions of the covariance structure of the underlying stationary process, with different observation lengths giving rise to a nested family of covariance matrices. These matrices are often compared using matrix distances, among which the affine-invariant Riemannian metric (AIRM) is particularly natural because it measures relative covariance differences and is invariant under nonsingular congruence transformations \cite{bhatia2009positive,pennec2006riemannian,moakher2005differential}. This raises a fundamental question: can the geometry of finite-dimensional covariance representations reveal an intrinsic geometric structure of the underlying infinite-dimensional process?
	
	A direct construction on infinite observation sequences is not straightforward. In particular, a nontrivial WSS process generally has infinite expected total energy and therefore does not define an ordinary finite-energy random element of \(\ell^2(\mathbb Z)\). Moreover, distinct ergodic stationary Gaussian models can induce mutually singular probability measures on the infinite observation space \cite{ibragimov2012gaussian,rozanov1963equivalence}. These facts make it important to distinguish the infinite-dimensional covariance structure from the regular finite-window statistical models used to extract information from observations. Rather than assigning a conventional likelihood to an infinite observation vector, we use finite-dimensional observation windows as compatible representations of a common underlying WSS covariance structure.
	
	The first part of this paper develops this construction for zero-mean proper complex Gaussian WSS processes. We formulate the covariance structure as a continuous positive semidefinite form on a space of discrete test functions. Its restriction to finite-dimensional observation subspaces produces the familiar Toeplitz covariance matrices, with the covariance matrix for a shorter window appearing as a leading principal submatrix of that for a longer window. This nested structure makes explicit that the finite-dimensional models are generated by one stationary covariance structure rather than being unrelated matrix models.
	
Using the Fisher information of each finite-window Gaussian likelihood, we first characterize its growth with the observation dimension. We show that the Fisher information grows linearly with the window length and study its normalized limit, which induces an intrinsic information metric on the space of stationary covariance models. Under standard Toeplitz regularity assumptions, this metric takes the form of a bilinear metric on spectral perturbations. In the logarithmic coordinate \(u=\log S\), the limiting metric reduces to the standard \(L^2\) inner product. Consequently, the Gaussian spectral Fisher--Rao geometry is flat in log-spectral coordinates, with an explicit geodesic distance. We further establish an asymptotic consistency result: the normalized squared AIRM distance between finite-window Toeplitz covariance matrices converges to the squared spectral Fisher--Rao distance. Thus, the spectral geometry is not introduced independently of finite-dimensional covariance geometry, but emerges as its large-window limit.
	
	The second part asks whether this spectral description can be extended beyond Gaussian observations. We consider finite-dimensional proper complex elliptical observation models whose covariance matrices are induced by the same type of WSS spectral structure. Elliptical models allow radial distributions that differ from the Gaussian case while retaining a covariance-based description of second-order dependence. Their geometry, however, contains an additional radial contribution. To isolate the part associated with covariance structure, we decompose each covariance matrix into a global scale and a determinant-normalized shape component.
	
	This decomposition leads to a Fisher-orthogonal separation between shape directions and nuisance directions associated with global scale and radial parameters. The shape metric is proportional to the canonical affine-invariant metric on the determinant-one covariance manifold, with a coefficient determined by the radial score. We characterize the asymptotic recovery of the Gaussian coefficient through the location Fisher information \(J_{N,\theta}\) of the log-radius. In particular, the coefficient converges to its Gaussian value when \(J_{N,\theta}=o(N^2)\). Common-scale Gaussian mixtures satisfy this condition. Here the same random scale is shared by the entire underlying Gaussian process; this detail ensures that the finite-dimensional distributions are marginals of one process and are therefore projectively consistent.
	
	Combining this radial result with Toeplitz asymptotics yields a limiting geometry on spectral shapes. In centered log-spectral coordinates, the limiting metric is the \(L^2\) inner product on zero-mean functions. Its geodesic distance is invariant under independent positive rescalings of the two spectra and agrees with the asymptotic normalized AIRM distance between their determinant-normalized covariance matrices. The resulting construction therefore separates global power from spectral shape while retaining a direct connection to finite-window covariance geometry.
	
	The contributions of this work can be summarized as follows:
	\begin{itemize}
		\item We formulate the covariance structure of WSS processes as a positive semidefinite bilinear form on a discrete test-function space, and establish its connection to compatible finite-dimensional Toeplitz covariance representations.
		\item We construct a spectral Fisher information rate for proper complex Gaussian WSS models and identify its flat \(L^2\) geometry in log-spectral coordinates. We further establish the asymptotic consistency between the spectral Fisher--Rao distance and the normalized AIRM distance between finite-dimensional Toeplitz covariance matrices.
		\item We extend the analysis to elliptical observation models by deriving a Fisher-orthogonal scale--shape decomposition, identifying a radial Fisher-information condition for recovery of the Gaussian shape coefficient, and showing that common-scale Gaussian mixtures provide a projectively consistent non-Gaussian class satisfying this condition.
	\end{itemize}
	
	The remainder of the paper is organized as follows. Section~2 introduces the weak covariance formulation of WSS processes and its spectral representation. Section~3 develops the Gaussian spectral information geometry and establishes its consistency with finite-window AIRM. Section~4 extends the construction to elliptical observations and derives the asymptotic spectral-shape geometry. Section~5 presents numerical experiments illustrating the principal theoretical results, and Section~6 concludes the paper.

	\section{Covariance Geometry and WSS Processes}
	\label{sec:covariance_geometry}
	
	Information geometry provides a geometric framework for
	characterizing and evaluating the distinction between statistical models.
	By equipping a statistical manifold with a Riemannian metric, the
	dissimilarity between probability distributions can be quantified
	through their geodesic distance \cite{amari2000methods,ay2017parametrized,nielsen2020elementary}.
	
	In statistical signal processing, many statistical models are characterized
	by their second-order statistics. For zero-mean real Gaussian distributions
	and zero-mean proper complex Gaussian distributions, the covariance matrix
	suffices to characterize the corresponding distribution \cite{neeser1993proper,picinbono1996second,schreier2010statistical}. The geometry of
	positive-definite covariance matrices therefore provides a natural means
	of comparing signals and noise conditions.
	
	The most commonly used metric is the affine-invariant Riemannian metric
	(AIRM), which measures relative differences between Hermitian positive-definite
	(HPD) matrices while remaining invariant under nonsingular congruence
	transformations \cite{bhatia2009positive,pennec2006riemannian,moakher2005differential}. For two HPD matrices
	$\mathbf{R}_N,\mathbf{S}_N\in\mathcal{P}_N$, the AIRM distance is defined as
	\begin{equation}
		d_{\mathrm{AIRM}}
		(\mathbf{R}_N,\mathbf{S}_N)
		=
		\left\|
		\log\left(
		\mathbf{R}_N^{-1/2}
		\mathbf{S}_N
		\mathbf{R}_N^{-1/2}
		\right)
		\right\|_{\mathrm F},
		\label{eq:airm_definition}
	\end{equation}
	where $\mathcal{P}_N$ denotes the manifold of $N\times N$
	Hermitian positive-definite matrices, $\log(\cdot)$ is the matrix
	logarithm, and $\|\cdot\|_{\mathrm F}$ is the Frobenius norm.
	
	A key property of the AIRM is its invariance under nonsingular
	congruence transformations:
	\begin{equation}
		d_{\mathrm{AIRM}}(\mathbf{R}_N,\mathbf{S}_N)
		=
		d_{\mathrm{AIRM}}
		\left(
		\mathbf{A}\mathbf{R}_N\mathbf{A}^{\mathrm H},
		\mathbf{A}\mathbf{S}_N\mathbf{A}^{\mathrm H}
		\right),
		\label{eq:airm_congruence}
	\end{equation}
	for any nonsingular matrix $\mathbf{A}$.
	
	\subsection{Weak Formulation of WSS Covariance Structures}
	
	Let $z[n]$ be a zero-mean, complex-valued random process
	with finite second-order moment at each time instant. Its
	autocorrelation function is defined as
	\[
	r[m,n] = \mathbb{E}\left[z[m]z[n]^*\right].
	\]
	
	We say that the process is Wide-Sense Stationary (WSS) if
	its autocorrelation function is invariant under a common
	time shift, i.e.,
	\[
	r[m+k,n+k]=r[m,n],\qquad \forall m,n,k\in\mathbb{Z}.
	\]
	In particular, setting $m=n$ shows that
	\[
	\mathbb{E}[|z[n]|^2]=\sigma^2
	\]
	is constant for all $n$. Consequently, any nontrivial
	infinite-length WSS process has infinite expected total
	energy:
	\[
	\mathbb{E}\left[\sum_{n\in\mathbb{Z}}|z[n]|^2\right]
	=
	\sum_{n\in\mathbb{Z}}\mathbb{E}[|z[n]|^2]
	=
	\infty,
	\]
	provided that $\sigma^2>0$. Thus, a nontrivial WSS process
	cannot, in general, be treated as a finite-energy random
	sequence in $\ell^2(\mathbb{Z})$.
	
	To formulate its covariance structure without restricting
	the process to a finite observation window, we first introduce
	the discrete Schwartz space
	\[
	\mathcal{S}(\mathbb{Z})
	=
	\left\{
	\varphi:\mathbb{Z}\rightarrow\mathbb{C}
	\;\middle|\;
	\sup_{n\in\mathbb{Z}}
	(1+|n|)^p|\varphi[n]|<\infty,
	\quad \forall p\geq 0
	\right\}.
	\]
	Every element of $\mathcal{S}(\mathbb{Z})$ belongs to both
	$\ell^1(\mathbb{Z})$ and $\ell^2(\mathbb{Z})$. Consequently, for any
	$\varphi\in\mathcal{S}(\mathbb{Z})$, the random series
	\[
	Z(\varphi)
	=
	\sum_{n\in\mathbb{Z}}z[n]\varphi[n]
	\]
	is well-defined. Indeed, by the Cauchy--Schwarz inequality,
	\[
	\begin{aligned}
		\mathbb{E}\left[
		\sum_{n\in\mathbb{Z}}|z[n]\varphi[n]|
		\right]
		&=
		\sum_{n\in\mathbb{Z}}
		\mathbb{E}[|z[n]|]|\varphi[n]|\\
		&\leq
		\sqrt{r[0]}\|\varphi\|_{\ell^1}
		<\infty.
	\end{aligned}
	\]
	Thus, the series converges absolutely almost surely. Moreover,
	$Z(\varphi)$ is square-integrable. In particular, by the
	Minkowski inequality,
	\[
	\mathbb{E}\left[|Z(\varphi)|^2\right]^{1/2}
	\leq
	\sum_{n\in\mathbb{Z}}
	\mathbb{E}[|z[n]|^2]^{1/2}|\varphi[n]|
	=
	\sqrt{r[0]}\|\varphi\|_{\ell^1},
	\]
	and hence
	\[
	\mathbb{E}\left[|Z(\varphi)|^2\right]
	\leq
	r[0]\|\varphi\|_{\ell^1}^2.
	\]
	
	We can therefore define the covariance form on $\mathcal{S}(\mathbb{Z})$ by
	\[
	\mathcal{B}(\varphi,\psi)
	=
	\mathbb{E}\left[
	Z(\varphi)^*Z(\psi)
	\right],
	\qquad
	\varphi,\psi\in\mathcal{S}(\mathbb{Z}).
	\]
	By construction, $\mathcal{B}$ is conjugate-linear in its first argument,
	linear in its second argument, and Hermitian:
	\[
	\mathcal{B}(\varphi,\psi)
	=
	\mathcal{B}(\psi,\varphi)^*.
	\]
	Furthermore, it is positive semidefinite, since
	\[
	\mathcal{B}(\varphi,\varphi)
	=
	\mathbb{E}\left[|Z(\varphi)|^2\right]
	\geq 0.
	\]
	The Cauchy--Schwarz inequality also gives
	\[
	|\mathcal{B}(\varphi,\psi)|
	\leq
	r[0]\|\varphi\|_{\ell^1}
	\|\psi\|_{\ell^1}.
	\]
	Since the $\ell^1$-norm is continuous with respect to the topology of
	$\mathcal{S}(\mathbb{Z})$, the covariance form is jointly continuous
	on this test-function space.
	
	This weak formulation allows the covariance structure of a WSS process
	to be represented independently of any particular finite observation
	window. In particular, let $\delta_n$ denote the discrete impulse at
	index $n$. By denoting
	\[
	z[n]=Z(\delta_n),
	\]
	the covariance form now satisfies
	\[
	\mathcal{B}(\delta_n,\delta_m)
	=
	\mathbb{E}[z[n]^*z[m]]
	=
	r[m-n].
	\]
	Hence, finite-dimensional covariance matrices can be recovered by
	restricting the same covariance form to finite-dimensional subspaces
	of test functions.
	
	More precisely, for
	\[
	\mathcal{H}_N
	=
	\operatorname{span}
	\{\delta_0,\delta_1,\ldots,\delta_{N-1}\},
	\]
	the Gram matrix of $\mathcal{B}$ with respect to the canonical basis is
	\[
	\mathbf{R}_N
	=
	\left[r[m-n]\right]_{n,m=0}^{N-1}.
	\]
	
	which is exactly the covariance matrix of dimension N.
	
\begin{figure}[h]
		\centering
		\includegraphics[width=0.6\linewidth]
		{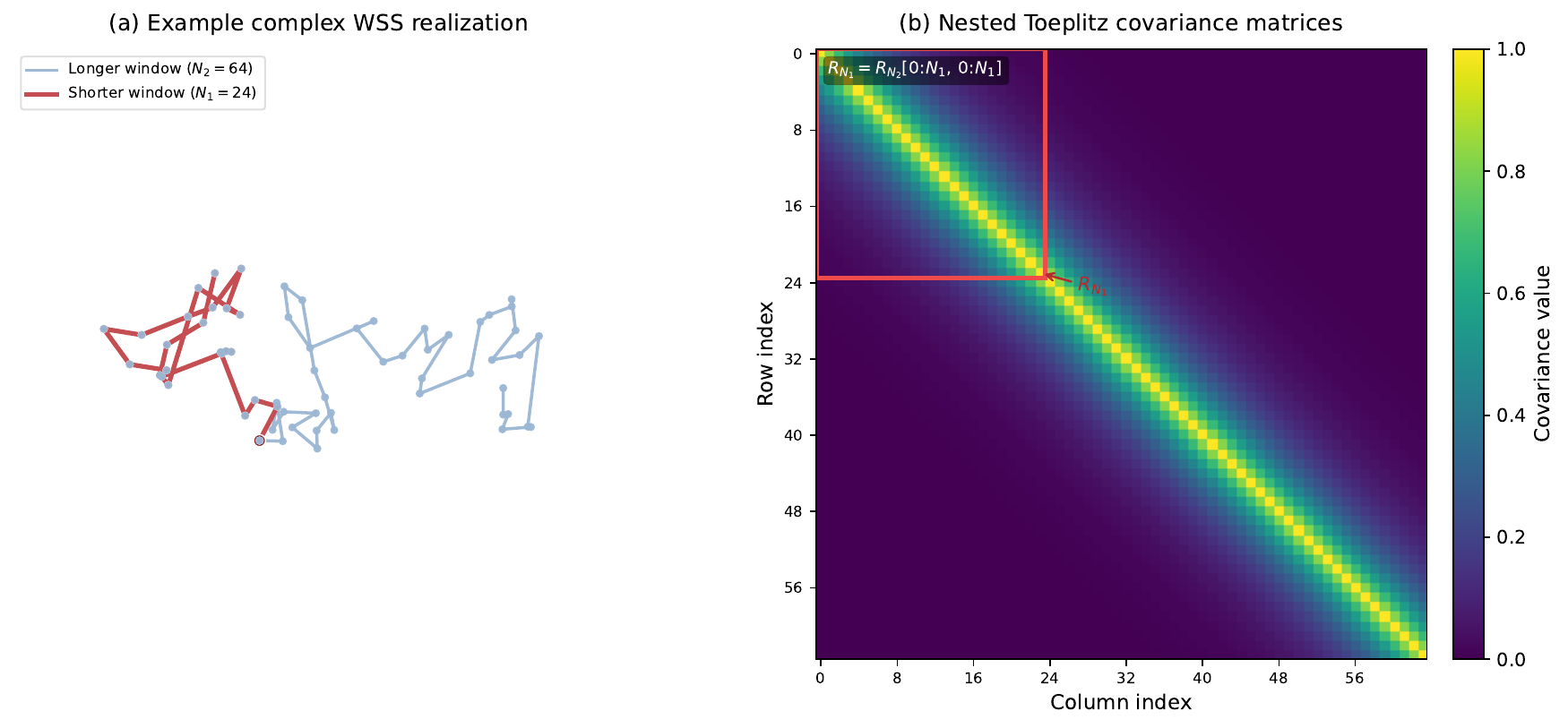}
		\caption{
			Illustration of the relationship between an infinite-dimensional
			WSS process and its finite-dimensional covariance representations.
			(a) An example realization of a proper complex Gaussian WSS process,
			with two nested observation windows of lengths $N_1$ and $N_2$.
			The three-dimensional trajectory is plotted using time, real part,
			and imaginary part as coordinates.
			(b) The corresponding theoretical Toeplitz covariance matrix
			$R_{N_2}$, with its leading principal submatrix $R_{N_1}$
			highlighted.
		}
		\label{fig:wss_covariance}
	\end{figure}

	The finite-dimensional covariance matrices induced by the
	underlying WSS process are naturally compatible across
	different observation lengths. In particular, for
	$N_1<N_2$, the covariance matrix associated with the
	shorter observation window is the leading principal
	submatrix of that associated with the longer window:
	\begin{equation}
		R_{N_1}
		=
		R_{N_2}[0:N_1,\;0:N_1].
		\label{eq:covariance_nesting}
	\end{equation}
	This nesting property reflects the fact that all
	finite-dimensional covariance models are induced by
	the same underlying autocorrelation structure.
	Figure~\ref{fig:wss_covariance} illustrates this
	relationship.
	
	\subsection{Spectral Representation of the Covariance Form}
	
	The weak covariance formulation admits a natural representation in the
	frequency domain. Since the autocorrelation sequence of a WSS process
	is positive definite, the Herglotz theorem guarantees the existence of \cite{brockwell2009time,ibragimov2012gaussian}
	a finite nonnegative measure $\mu$ on the frequency torus
	$\mathbb{T}=[-\pi,\pi)$ such that
	\[
	r[k]
	=
	\int_{-\pi}^{\pi}e^{jk\omega}\,\mu(d\omega),
	\qquad k\in\mathbb{Z},
	\]
	where
	\[
	\mu(\mathbb{T})=r[0].
	\]
	This representation holds for general WSS processes and does not require
	the autocorrelation sequence to be absolutely summable. In particular,
	the spectral measure may contain both absolutely continuous and
	singular components.
	
	For a test function $\varphi\in\mathcal{S}(\mathbb{Z})$, define its
	Fourier transform by
	\[
	\widehat{\varphi}(\omega)
	=
	\sum_{n\in\mathbb{Z}}\varphi[n]e^{j\omega n}.
	\]
	Since $\varphi$ belongs to the discrete Schwartz space, its Fourier
	transform is a smooth periodic function. Substituting the spectral
	representation of $r[k]$ into the covariance form gives
	\[
	\begin{aligned}
		\mathcal{B}(\varphi,\psi)
		&=
		\sum_{n,m\in\mathbb{Z}}
		r[m-n]\varphi[n]^*\psi[m]\\
		&=
		\int_{-\pi}^{\pi}
		\left(
		\sum_{n\in\mathbb{Z}}\varphi[n]^*e^{-jn\omega}
		\right)
		\left(
		\sum_{m\in\mathbb{Z}}\psi[m]e^{jm\omega}
		\right)
		\mu(d\omega)\\
		&=
		\int_{-\pi}^{\pi}
		\widehat{\varphi}(\omega)^*
		\widehat{\psi}(\omega)\,\mu(d\omega).
	\end{aligned}
	\]
	Therefore, the covariance form is represented in the frequency domain as
	\[
	\boxed{
		\mathcal{B}(\varphi,\psi)
		=
		\int_{-\pi}^{\pi}
		\widehat{\varphi}(\omega)^*
		\widehat{\psi}(\omega)\,\mu(d\omega).
	}
	\]
	
	When the spectral measure is absolutely continuous with respect to
	Lebesgue measure, it admits a density $f(\omega)$, so that
	\[
	\mu(d\omega)=f(\omega)\,d\omega.
	\]
	Introducing the power spectral density
	\[
	S(\omega)=2\pi f(\omega),
	\]
	the covariance form becomes
	\[
	\mathcal{B}(\varphi,\psi)
	=
	\frac{1}{2\pi}
	\int_{-\pi}^{\pi}
	S(\omega)
	\widehat{\varphi}(\omega)^*
	\widehat{\psi}(\omega)\,d\omega.
	\]
	In this case, the frequency-domain covariance structure is characterized
	by multiplication by the spectral density $S(\omega)$. More precisely,
	under the Fourier convention adopted above, the covariance operator
	associated with the form acts in the frequency domain as multiplication
	by $S(\omega)$.
	
	In this work, we focus on the absolutely continuous case
	with strictly positive and bounded spectral densities.
	This assumption covers a broad class of practical
	signal models in which a non-vanishing noise floor keeps
	the spectrum strictly positive. It also ensures that the
	natural logarithm of the spectral density is well-defined
	and bounded.
	
	\section{Spectral-Domain Information Geometry}
	\label{sec:gaussian_model}
	
	In the previous section, we established a frequency-domain representation
	of the covariance structure of an infinite-length WSS process.
	In particular, when the spectral measure is absolutely continuous,
	the covariance structure is completely characterized by its power
	spectral density (PSD).
	
	The goal of this section is to develop an information-geometric framework
	directly on the space of spectral densities. We focus on zero-mean
	proper complex Gaussian WSS processes and follow the Fisher--Rao
	geometric framework to characterize the statistical distinguishability
	between stationary models.
	
	The first natural attempt is to work directly with the probability measures of infinite-length stationary processes. For stationary Gaussian processes, the absence of spectral atoms ensures ergodicity. Consequently, distinct stationary Gaussian models with absolutely continuous spectral densities induce mutually singular probability measures on the infinite observation space. This leads to an all-or-nothing statistical separation between distinct models, preventing the conventional likelihood-based Fisher--Rao construction from providing a nontrivial local characterization of their distinguishability.
	
	A natural alternative is to begin with finite observation windows, for which the corresponding Gaussian statistical models admit well-defined likelihood functions and Fisher information matrices. As the observation length increases, the Fisher information generally grows, reflecting the accumulation of statistical information over time. We therefore investigate its asymptotic growth and identify an appropriate normalization to remove the dependence on the observation length, with the aim of obtaining a finite and nontrivial limiting characterization of the statistical information.
	
	Based on this construction, we derive a Riemannian metric directly on the space of spectral densities. The resulting geometry characterizes the local statistical distinguishability of stationary Gaussian models in terms of their spectral representations. We further derive its geodesic distance and establish its asymptotic consistency with the affine-invariant Riemannian distance (AIRM) between finite-dimensional covariance matrices.
	
	\subsection{Gaussian WSS Statistical Model}
	
	Let \(z[n]\) be a zero-mean proper complex Gaussian WSS process defined for all \(n\in\mathbb Z\). We assume that its spectral measure is absolutely continuous with respect to Lebesgue measure, such that its covariance function admits the representation
	\[
	r[k]
	=
	\frac{1}{2\pi}
	\int_{-\pi}^{\pi}
	S(\omega)e^{jk\omega}\,d\omega,
	\]
	where \(S(\omega)\) denotes the power spectral density (PSD).
	
	Throughout this section, we consider strictly positive and bounded spectral densities belonging to
	\[
	\mathcal P
	=
	\left\{
	S\in C(\mathbb T;\mathbb R)
	\;\middle|\;
	0<m_S\leq S(\omega)\leq M_S<\infty
	\ \text{a.e.}
	\right\},
	\]
	where \(m_S\) and \(M_S\) may depend on \(S\).
	
	Since a zero-mean proper complex Gaussian process is completely determined by its covariance function, each \(S\in\mathcal P\) uniquely specifies an infinite-length stationary Gaussian model, denoted by \(P_S\). Thus, \(\mathcal P\) serves as the parameter space of the statistical models considered in this work.
	
	\subsection{Fisher Information Rate on the Spectral Density Space}
	Consider an observation window \(\{z[0],z[1],\ldots,z[N-1]\}\) of the process. To maintain the covariance convention adopted in the preceding section, define the conjugated observation vector
	\[
	\mathbf{y}_N
	=
	[z[0]^*,z[1]^*,\ldots,z[N-1]^*]^{\mathsf T}.
	\]
	
	Although each spectral density specifies an infinite-length statistical model, its finite-dimensional observation vectors provide regular likelihood-based representations of that model. In particular, \(\mathbf y_N\) is a zero-mean proper complex Gaussian random vector whose covariance matrix is
	\[
	\mathbf R_N(S)
	=
	\left[
	\frac{1}{2\pi}
	\int_{-\pi}^{\pi}
	S(\omega)e^{j(m-n)\omega}\,d\omega
	\right]_{n,m=0}^{N-1}.
	\]
	This matrix is induced by the spectral density \(S\) and the chosen observation window. In particular, changing \(N\) does not change the underlying statistical model.

	The probability density of $\mathbf{y}_N$ is
	\begin{equation}
		p_S(\mathbf{y}_N)
		=
		\frac{1}{\pi^N\det\mathbf{R}_N(S)}
		\exp\left(
		-\mathbf{y}_N^{\mathrm H}
		\mathbf{R}_N(S)^{-1}
		\mathbf{y}_N
		\right).
		\label{eq:complex_gaussian_likelihood}
	\end{equation}
	
	Consider a differentiable perturbation of the spectral density,
	\begin{equation}
		S_t=S+tH,
	\end{equation}
	where $H$ is a real-valued continuous function and $t$ is sufficiently
	small such that $S_t$ remains strictly positive. The corresponding
	covariance perturbation is
	\begin{equation}
		\mathbf{R}_{N,H}
		=
		\left.
		\frac{d}{dt}\mathbf{R}_N(S+tH)
		\right|_{t=0}
		=
		\mathbf{R}_N(H).
	\end{equation}
	
	Differentiating the log-likelihood in the direction $H$ gives
	the score function
	\begin{equation}
		\begin{aligned}
			\ell_{N,H}(\mathbf{y}_N)
			&=
			\left.
			\frac{d}{dt}
			\log p_{S+tH}(\mathbf{y}_N)
			\right|_{t=0}\\
			&=
			-\operatorname{tr}
			\left(
			\mathbf{R}_N^{-1}\mathbf{R}_{N,H}
			\right)
			+
			\mathbf{y}_N^{\mathrm H}
			\mathbf{R}_N^{-1}
			\mathbf{R}_{N,H}
			\mathbf{R}_N^{-1}
			\mathbf{y}_N,
		\end{aligned}
		\label{eq:gaussian_score}
	\end{equation}
	where $\mathbf{R}_N=\mathbf{R}_N(S)$.
	
	For two perturbations $H$ and $K$, the Fisher information is
	the expected product of their score functions. Using the
	fourth-order moment identity for proper complex Gaussian random
	vectors, we obtain
	\begin{equation}
		\begin{aligned}
			g_{N,S}(H,K)
			&=
			\mathbb{E}_S
			\left[
			\ell_{N,H}(\mathbf{y}_N)
			\ell_{N,K}(\mathbf{y}_N)
			\right]\\
			&=
			\operatorname{tr}
			\left(
			\mathbf{R}_N^{-1}
			\mathbf{R}_N(H)
			\mathbf{R}_N^{-1}
			\mathbf{R}_N(K)
			\right).
		\end{aligned}
		\label{eq:finite_fisher_metric}
	\end{equation}
	This is the Fisher information of the finite-window observation
	model. It also coincides with the local form of the AIRM on the
	corresponding positive-definite covariance matrices.
	
	The Fisher information in \ref{eq:finite_fisher_metric} exhibits two distinct contributions as the observation length \(N\) varies: a leading-order term that grows linearly with \(N\), reflecting the accumulation of statistical information, while finite-window effects are captured by the asymptotic approximation error. Our goal is to obtain a geometric characterization that is independent of the observation length. To this end, we investigate the asymptotic behaviour of the finite-window Fisher information as \(N\) increases.
	
	Define
	\begin{equation}
		g_{N,S}^{\mathrm{norm}}(H,K)
		=
		\frac{1}{N}g_{N,S}(H,K).
		\label{eq:normalized_fisher_metric}
	\end{equation}
	
	We next look for the limit of this expression:
	
	For a continuous function $a$ on $\mathbb{T}$, denote by
	$\mathbf{T}_N(a)$ the Toeplitz matrix
	\begin{equation}
		[\mathbf{T}_N(a)]_{n,m}
		=
		\frac{1}{2\pi}
		\int_{-\pi}^{\pi}
		a(\omega)e^{j(m-n)\omega}\,d\omega.
	\end{equation}
	Then
	\begin{equation}
		\mathbf{R}_N(S)=\mathbf{T}_N(S),\qquad
		\mathbf{R}_N(H)=\mathbf{T}_N(H).
	\end{equation}
	
	Assume that $S$ is continuous and uniformly positive, and that
	$H$ and $K$ are continuous. The asymptotic functional calculus
	for Toeplitz matrices gives
	\begin{equation}
		\left\|
		\mathbf{T}_N(S)^{-1}
		-
		\mathbf{T}_N(S^{-1})
		\right\|_{2,N}
		\longrightarrow 0,
		\label{eq:toeplitz_inverse_approx}
	\end{equation}
	where
	\begin{equation}
		\|\mathbf{X}_N\|_{2,N}
		=
		\frac{1}{\sqrt{N}}\|\mathbf{X}_N\|_{\mathrm F}.
	\end{equation}
	Furthermore, the asymptotic product property of Toeplitz matrices
	implies
	\begin{equation}
		\begin{aligned}
			&\left\|
			\mathbf{T}_N(S)^{-1}
			\mathbf{T}_N(H)
			\mathbf{T}_N(S)^{-1}
			\mathbf{T}_N(K)
			-
			\mathbf{T}_N\left(\frac{HK}{S^2}\right)
			\right\|_{2,N}\\
			&\hspace{4cm}\longrightarrow 0.
		\end{aligned}
		\label{eq:toeplitz_fisher_product}
	\end{equation}
	
	The uniform boundedness of the operator norms of the matrices
	involved allows the normalized trace to be evaluated through
	this approximation. Consequently,
	\begin{equation}
		\begin{aligned}
			&\lim_{N\to\infty}
			\frac{1}{N}
			\operatorname{tr}
			\left[
			\mathbf{R}_N(S)^{-1}
			\mathbf{R}_N(H)
			\mathbf{R}_N(S)^{-1}
			\mathbf{R}_N(K)
			\right]\\
			&=
			\lim_{N\to\infty}
			\frac{1}{N}
			\operatorname{tr}
			\mathbf{T}_N\left(\frac{HK}{S^2}\right)\\
			&=
			\frac{1}{2\pi}
			\int_{-\pi}^{\pi}
			\frac{H(\omega)K(\omega)}
			{S^2(\omega)}\,d\omega.
		\end{aligned}
		\label{eq:fisher_rate_limit}
	\end{equation}
	The last equality follows directly from the diagonal entries
	of the Toeplitz matrix.
	
	Thus, the normalized Fisher information converges to a bilinear
	form that depends only on the spectral density and its perturbations,
	and not on the observation length. We define the resulting
	spectral Fisher information rate as
	\begin{equation}
		\boxed{
			g_S(H,K)
			=
			\frac{1}{2\pi}
			\int_{-\pi}^{\pi}
			\frac{H(\omega)K(\omega)}
			{S^2(\omega)}\,d\omega.
		}
		\label{eq:spectral_fisher_metric}
	\end{equation}
	
	Equation \eqref{eq:spectral_fisher_metric} defines the local information geometry of stationary Gaussian models. This geometry is intrinsic to the spectral density and independent of the observation length.
	
	\subsection{Geodesic Distance and Log-Spectral Geometry}
	
	We now investigate the geometry induced by the Fisher information
	rate in \eqref{eq:spectral_fisher_metric}. The parameter space
	$\mathcal{P}$ is an open subset of the real Banach space
	$C(\mathbb{T};\mathbb{R})$. Its tangent space at any
	$S\in\mathcal{P}$ can therefore be identified with
	\begin{equation}
		T_S\mathcal{P}
		=
		C(\mathbb{T};\mathbb{R}).
	\end{equation}
	
	For any nonzero $H\in T_S\mathcal{P}$, the metric satisfies
	\begin{equation}
		g_S(H,H)
		=
		\frac{1}{2\pi}
		\int_{-\pi}^{\pi}
		\frac{H^2(\omega)}{S^2(\omega)}\,d\omega
		>0.
	\end{equation}
	Moreover, since $S$ is bounded and bounded away from zero,
	the norm induced by this metric is equivalent to the $L^2$
	norm on tangent vectors. It is therefore weaker than the ambient $C(\mathbb T)$ topology, and $g_S$ defines a weak Riemannian metric
	on $\mathcal{P}$.
	
	A particularly useful property of this metric is revealed by
	the logarithmic coordinate transformation
	\begin{equation}
		u(\omega)=\log S(\omega).
		\label{eq:log_spectral_coordinate}
	\end{equation}
	The logarithmic map is a smooth bijection from $\mathcal{P}$
	onto $C(\mathbb{T};\mathbb{R})$, with inverse given by
	the pointwise exponential map.
	
	For a differentiable curve $S_t\in\mathcal{P}$, let
	\begin{equation}
		H_t=\partial_t S_t.
	\end{equation}
	The corresponding velocity in logarithmic coordinates is
	\begin{equation}
		\partial_t u_t(\omega)
		=
		\frac{H_t(\omega)}{S_t(\omega)}.
	\end{equation}
	Substituting this relation into the Fisher metric gives
	\begin{equation}
		\begin{aligned}
			g_{S_t}(H_t,H_t)
			&=
			\frac{1}{2\pi}
			\int_{-\pi}^{\pi}
			\frac{H_t^2(\omega)}{S_t^2(\omega)}\,d\omega\\
			&=
			\frac{1}{2\pi}
			\int_{-\pi}^{\pi}
			|\partial_t u_t(\omega)|^2\,d\omega.
		\end{aligned}
		\label{eq:log_spectral_metric}
	\end{equation}
	Hence, the logarithmic coordinate transformation maps the spectral
	Fisher geometry to the standard flat $L^2$ geometry.
	
	It follows that the shortest path between two points in logarithmic
	coordinates is the straight line connecting them. For
	$S_1,S_2\in\mathcal{P}$, define
	\begin{equation}
		u_t(\omega)
		=
		(1-t)\log S_1(\omega)
		+
		t\log S_2(\omega),
		\qquad 0\leq t\leq 1.
	\end{equation}
	Mapping this path back to the spectral density space gives
	the geodesic
	\begin{equation}
		\boxed{
			S_t(\omega)
			=
			S_1^{1-t}(\omega)S_2^t(\omega),
			\qquad 0\leq t\leq 1.
		}
		\label{eq:spectral_geodesic}
	\end{equation}
	
	The length of this geodesic is
	\begin{equation}
		\begin{aligned}
			L(S_t)
			&=
			\int_0^1
			\sqrt{g_{S_t}(\partial_tS_t,\partial_tS_t)}
			\,dt\\
			&=
			\left[
			\frac{1}{2\pi}
			\int_{-\pi}^{\pi}
			\left|
			\log S_2(\omega)-\log S_1(\omega)
			\right|^2d\omega
			\right]^{1/2}.
		\end{aligned}
	\end{equation}
	Since the logarithmic representation is an isometry into the
	flat $L^2$ space, this path has minimum length among all
	differentiable paths connecting $S_1$ and $S_2$.
	
	Therefore, the geodesic distance induced by the Fisher information
	rate is
	\begin{equation}
		\boxed{
			d_{\mathrm{FR}}^2(S_1,S_2)
			=
			\frac{1}{2\pi}
			\int_{-\pi}^{\pi}
			\left|
			\log S_2(\omega)-\log S_1(\omega)
			\right|^2d\omega.
		}
		\label{eq:spectral_fisher_distance}
	\end{equation}

	The distance \(d_{\mathrm{FR}}\) satisfies the metric properties. In particular,
	\begin{equation}
		d_{\mathrm{FR}}(S_1,S_2)=0
		\quad\Longleftrightarrow\quad
		S_1=S_2\quad\text{a.e.},
	\end{equation}
	since the logarithm is injective on the positive real numbers. Symmetry and the triangle inequality follow directly from the corresponding properties of the \(L^2\) norm in logarithmic coordinates. Therefore, \(d_{\mathrm{FR}}\) defines a metric on \(\mathcal{P}\).
	
\begin{figure}[h]
		\centering
		\includegraphics[width=0.95\linewidth]
		{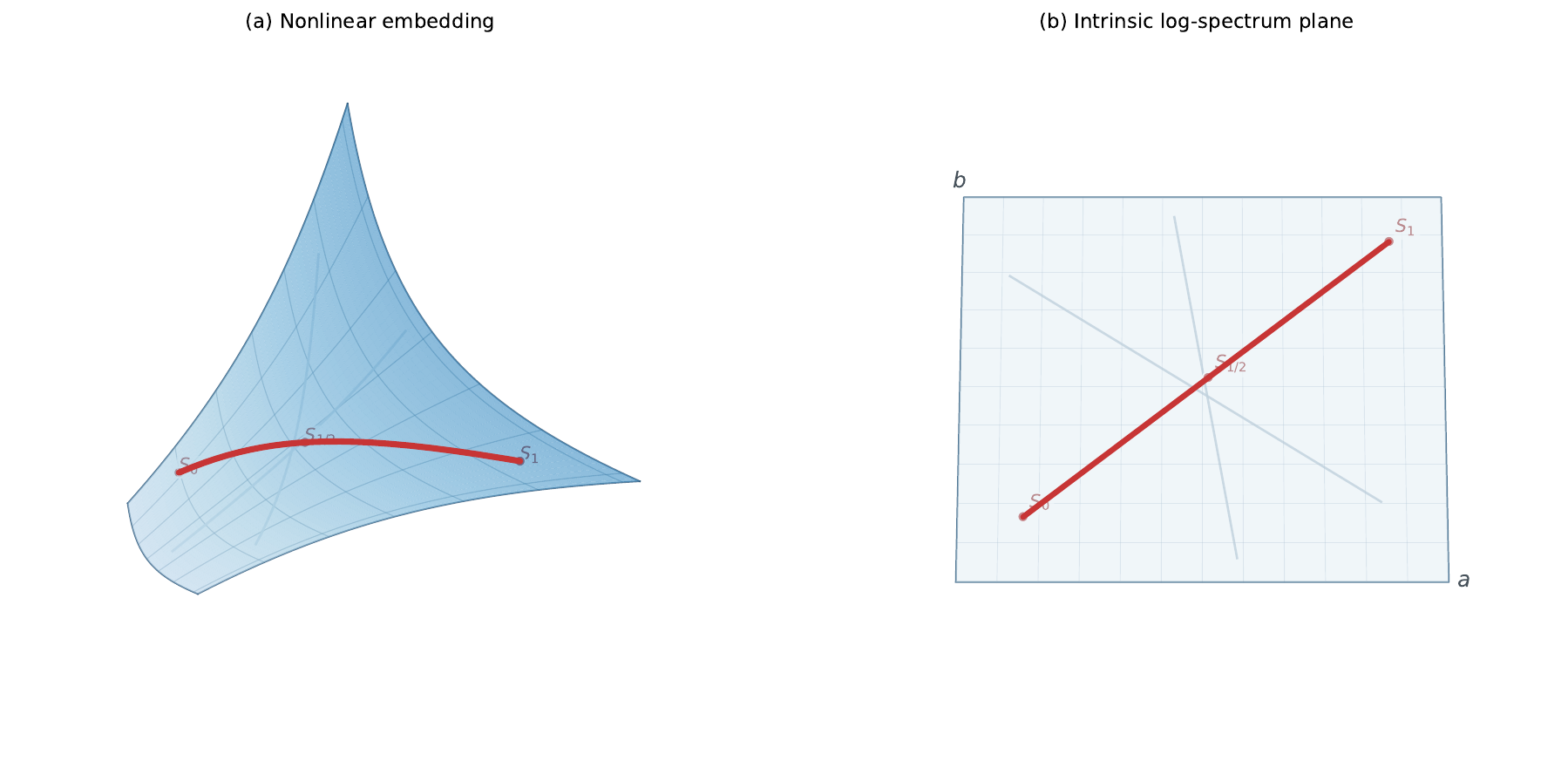}
		\caption{Geometric illustration of the Gaussian spectral
			Fisher--Rao geometry. (a) A nonlinear three-dimensional
			evaluation embedding of a two-parameter spectral family.
			The highlighted curve is the image of a Fisher--Rao geodesic
			under this embedding, which is not required to preserve the
			Fisher metric. (b) The same geodesic in the intrinsic
			parameter coordinates $(a,b)$, where the Fisher metric is
			Euclidean up to a constant factor and geodesics are straight
			lines. The auxiliary straight lines illustrate other geodesics
			in the same family.}
		\label{fig:spectral_geodesic}
	\end{figure}

For illustration, consider the two-parameter spectral family
\[
u(\omega;a,b)=a\cos\omega+b\cos 2\omega.
\]
Figure~\ref{fig:spectral_geodesic} illustrates the corresponding
geodesic in the nonlinear embedding and intrinsic log-spectrum
coordinates.
	
	\subsection{Asymptotic Consistency with Finite-Dimensional AIRM}
	
	We now show that the spectral Fisher--Rao
	distance is recovered asymptotically from the AIRM distances
	between these finite-window covariance representations.
	
	For $S\in\mathcal{P}$, define
	\begin{equation}
		\mathbf{R}_N(S)
		=
		\left[
		\frac{1}{2\pi}
		\int_{-\pi}^{\pi}
		S(\omega)e^{j(m-n)\omega}\,d\omega
		\right]_{n,m=0}^{N-1}.
		\label{eq:finite_toeplitz_covariance}
	\end{equation}
	For strictly positive real-valued $S$, this matrix is Hermitian
	positive definite. Importantly, the family
	$\{\mathbf{R}_N(S)\}_{N\geq1}$ consists of finite-window
	representations generated by the same spectral density, rather
	than independently specified covariance models.
	
	For two Hermitian positive-definite matrices
	$\mathbf{A},\mathbf{B}\in\mathcal{P}_N$, their affine-invariant
	Riemannian distance is
	\begin{equation}
		d_{\mathrm{AIRM}}(\mathbf{A},\mathbf{B})
		=
		\left\|
		\log\left(
		\mathbf{A}^{-1/2}
		\mathbf{B}
		\mathbf{A}^{-1/2}
		\right)
		\right\|_{\mathrm F}.
		\label{eq:finite_airm_definition}
	\end{equation}
	
	We now establish the asymptotic relationship between the
	finite-dimensional AIRM and the spectral Fisher--Rao distance.
	
	\begin{theorem}[Asymptotic consistency of finite-window AIRM]
		\label{thm:airm_consistency}
		Let $S_1,S_2$ be real-valued continuous functions on $\mathbb{T}$
		satisfying
		\begin{equation}
			0<m\leq S_i(\omega)\leq M<\infty,
			\qquad i=1,2,
		\end{equation}
		for all $\omega\in\mathbb{T}$. Let
		\begin{equation}
			\mathbf{R}_N=\mathbf{R}_N(S_1),
			\qquad
			\mathbf{S}_N=\mathbf{R}_N(S_2).
		\end{equation}
		Then,
		\begin{equation}
			\boxed{
				\lim_{N\to\infty}
				\frac{1}{N}
				d_{\mathrm{AIRM}}^2
				(\mathbf{R}_N,\mathbf{S}_N)
				=
				\frac{1}{2\pi}
				\int_{-\pi}^{\pi}
				\left|
				\log S_2(\omega)-\log S_1(\omega)
				\right|^2\,d\omega.
			}
			\label{eq:airm_spectral_convergence}
		\end{equation}
		Equivalently,
		\begin{equation}
			\lim_{N\to\infty}
			\frac{1}{\sqrt{N}}
			d_{\mathrm{AIRM}}(\mathbf{R}_N,\mathbf{S}_N)
			=
			d_{\mathrm{FR}}(S_1,S_2).
		\end{equation}
	\end{theorem}

	\begin{proof}
		
		For notational convenience, define the normalized Frobenius norm
		\begin{equation}
			\|\mathbf{X}_N\|_{2,N}
			=
			\frac{1}{\sqrt{N}}\|\mathbf{X}_N\|_{\mathrm F}.
		\end{equation}
		Let
		\begin{equation}
			\mathbf{A}_N=\mathbf{R}_N(S_1),
			\qquad
			\mathbf{B}_N=\mathbf{R}_N(S_2),
		\end{equation}
		and introduce
		\begin{equation}
			\mathbf{H}_N
			=
			\mathbf{A}_N^{-1/2}
			\mathbf{B}_N
			\mathbf{A}_N^{-1/2}.
		\end{equation}
		By the definition of the AIRM,
		\begin{equation}
			d_{\mathrm{AIRM}}^2(\mathbf{A}_N,\mathbf{B}_N)
			=
			\|\log\mathbf{H}_N\|_{\mathrm F}^2.
		\end{equation}
		
		We will show that the normalized logarithmic matrix
		$\log\mathbf{H}_N$ can be asymptotically approximated by the
		logarithm of a Toeplitz matrix generated by the pointwise spectral
		ratio $S_2/S_1$.
		
		\textit{Step 1: Approximation of the relative covariance matrix
			by a Toeplitz matrix.}
		
		For a continuous function $c$ on $\mathbb{T}$, denote by
		$\mathbf{T}_N(c)$ the Toeplitz matrix with entries
		\begin{equation}
			[\mathbf{T}_N(c)]_{n,m}
			=
			\frac{1}{2\pi}
			\int_{-\pi}^{\pi}
			c(\omega)e^{j(m-n)\omega}\,d\omega.
		\end{equation}
		Thus,
		\begin{equation}
			\mathbf{A}_N=\mathbf{T}_N(S_1),
			\qquad
			\mathbf{B}_N=\mathbf{T}_N(S_2).
		\end{equation}
		
		Define
		\begin{equation}
			\mathbf{C}_N
			=
			\mathbf{T}_N(S_1^{-1/2}).
		\end{equation}
		The asymptotic functional calculus for Toeplitz matrices gives \cite{gray2006toeplitz,bottcher2012introduction,grenander2026toeplitz}.
		\begin{equation}
			\left\|
			\mathbf{A}_N^{-1/2}-\mathbf{C}_N
			\right\|_{2,N}
			\longrightarrow 0.
		\end{equation}
		Furthermore, the operator norms satisfy
		\begin{equation}
			\|\mathbf{A}_N^{-1/2}\|_{\mathrm{op}}
			\leq m^{-1/2},
			\qquad
			\|\mathbf{C}_N\|_{\mathrm{op}}
			\leq m^{-1/2},
		\end{equation}
		and
		\begin{equation}
			\|\mathbf{B}_N\|_{\mathrm{op}}\leq M.
		\end{equation}
		
		Using the decomposition
		\begin{equation}
			\begin{aligned}
				\mathbf{H}_N-\mathbf{C}_N\mathbf{B}_N\mathbf{C}_N
				={}&
				(\mathbf{A}_N^{-1/2}-\mathbf{C}_N)
				\mathbf{B}_N\mathbf{A}_N^{-1/2}\\
				&+
				\mathbf{C}_N\mathbf{B}_N
				(\mathbf{A}_N^{-1/2}-\mathbf{C}_N),
			\end{aligned}
		\end{equation}
		together with the matrix norm inequality
		\begin{equation}
			\|\mathbf{X}\mathbf{Y}\mathbf{Z}\|_{\mathrm F}
			\leq
			\|\mathbf{X}\|_{\mathrm F}
			\|\mathbf{Y}\|_{\mathrm{op}}
			\|\mathbf{Z}\|_{\mathrm{op}},
		\end{equation}
		we obtain
		\begin{equation}
			\left\|
			\mathbf{H}_N-\mathbf{C}_N\mathbf{B}_N\mathbf{C}_N
			\right\|_{2,N}
			\longrightarrow 0.
		\end{equation}
		
		Next, the asymptotic product property of Toeplitz matrices gives \cite{gray2006toeplitz,bottcher2012introduction}.
		\begin{equation}
			\left\|
			\mathbf{T}_N(S_1^{-1/2})
			\mathbf{T}_N(S_2)
			-
			\mathbf{T}_N(S_1^{-1/2}S_2)
			\right\|_{2,N}
			\longrightarrow 0.
		\end{equation}
		Applying the same property once more, we obtain
		\begin{equation}
			\left\|
			\mathbf{T}_N(S_1^{-1/2}S_2)
			\mathbf{T}_N(S_1^{-1/2})
			-
			\mathbf{T}_N(S_2/S_1)
			\right\|_{2,N}
			\longrightarrow 0.
		\end{equation}
		Since the operator norms of the Toeplitz factors are uniformly
		bounded, these two relations imply
		\begin{equation}
			\left\|
			\mathbf{C}_N\mathbf{B}_N\mathbf{C}_N
			-
			\mathbf{T}_N(S_2/S_1)
			\right\|_{2,N}
			\longrightarrow 0.
		\end{equation}
		
		Consequently, defining
		\begin{equation}
			\mathbf{D}_N=\mathbf{T}_N(S_2/S_1),
		\end{equation}
		we have
		\begin{equation}
			\boxed{
				\|\mathbf{H}_N-\mathbf{D}_N\|_{2,N}
				\longrightarrow 0.
			}
			\label{eq:relative_toeplitz_approx}
		\end{equation}
		
		\textit{Step 2: Approximation of the matrix logarithms.}
		
		Since the symbols satisfy
		\begin{equation}
			m\leq S_i(\omega)\leq M,
		\end{equation}
		the eigenvalues of $\mathbf{A}_N$ and $\mathbf{B}_N$ are
		contained in $[m,M]$. It follows that
		\begin{equation}
			\frac{m}{M}\mathbf{I}_N
			\preceq
			\mathbf{H}_N
			\preceq
			\frac{M}{m}\mathbf{I}_N.
		\end{equation}
		Similarly, since
		\begin{equation}
			\frac{m}{M}
			\leq
			\frac{S_2(\omega)}{S_1(\omega)}
			\leq
			\frac{M}{m},
		\end{equation}
		the eigenvalues of $\mathbf{D}_N$ satisfy
		\begin{equation}
			\frac{m}{M}\mathbf{I}_N
			\preceq
			\mathbf{D}_N
			\preceq
			\frac{M}{m}\mathbf{I}_N.
		\end{equation}
		
		For two Hermitian positive-definite matrices $\mathbf{X}$ and
		$\mathbf{Y}$ satisfying
		\begin{equation}
			\mathbf{X},\mathbf{Y}\succeq\alpha\mathbf{I},
		\end{equation}
		the matrix logarithm obeys the Frobenius norm estimate
		\begin{equation}
			\|\log\mathbf{X}-\log\mathbf{Y}\|_{\mathrm F}
			\leq
			\frac{1}{\alpha}
			\|\mathbf{X}-\mathbf{Y}\|_{\mathrm F}.
		\end{equation}
		Applying it with $\alpha=m/M$, we obtain
		\begin{equation}
			\begin{aligned}
				\|\log\mathbf{H}_N-\log\mathbf{D}_N\|_{2,N}
				&\leq
				\frac{M}{m}
				\|\mathbf{H}_N-\mathbf{D}_N\|_{2,N}\\
				&\longrightarrow 0.
			\end{aligned}
		\end{equation}
		
		Furthermore, the uniform spectral bounds imply that
		\begin{equation}
			\|\log\mathbf{H}_N\|_{2,N}
			+
			\|\log\mathbf{D}_N\|_{2,N}
		\end{equation}
		is uniformly bounded. Therefore, using
		\begin{equation}
			\begin{aligned}
				\big|
				\|\mathbf{X}\|_{\mathrm F}^2
				-
				\|\mathbf{Y}\|_{\mathrm F}^2
				\big|
				\leq
				\|\mathbf{X}-\mathbf{Y}\|_{\mathrm F}
				\big(
				\|\mathbf{X}\|_{\mathrm F}
				+
				\|\mathbf{Y}\|_{\mathrm F}
				\big),
			\end{aligned}
		\end{equation}
		we conclude that
		\begin{equation}
			\frac{1}{N}
			\left|
			\|\log\mathbf{H}_N\|_{\mathrm F}^2
			-
			\|\log\mathbf{D}_N\|_{\mathrm F}^2
			\right|
			\longrightarrow 0.
			\label{eq:log_frobenius_convergence}
		\end{equation}
		
		\textit{Step 3: Application of the Szeg\H{o} eigenvalue
			distribution theorem.}
		
		The matrix $\mathbf{D}_N$ is a Hermitian Toeplitz matrix generated
		by the continuous positive function
		\begin{equation}
			c(\omega)=\frac{S_2(\omega)}{S_1(\omega)}.
		\end{equation}
		By the Szeg\H{o} eigenvalue distribution theorem \cite{gray2006toeplitz,grenander2026toeplitz}, for any continuous
		function $F$ on the spectral interval containing the eigenvalues
		of $\mathbf{D}_N$,
		\begin{equation}
			\lim_{N\to\infty}
			\frac{1}{N}\operatorname{tr}F(\mathbf{D}_N)
			=
			\frac{1}{2\pi}
			\int_{-\pi}^{\pi}F(c(\omega))\,d\omega.
		\end{equation}
		Taking
		\begin{equation}
			F(x)=(\log x)^2,
		\end{equation}
		which is continuous on $[m/M,M/m]$, gives
		\begin{equation}
			\lim_{N\to\infty}
			\frac{1}{N}
			\operatorname{tr}
			\left[(\log\mathbf{D}_N)^2\right]
			=
			\frac{1}{2\pi}
			\int_{-\pi}^{\pi}
			\left|\log\frac{S_2(\omega)}{S_1(\omega)}\right|^2
			\,d\omega.
		\end{equation}
		Since $\mathbf{D}_N$ is Hermitian positive definite,
		\begin{equation}
			\|\log\mathbf{D}_N\|_{\mathrm F}^2
			=
			\operatorname{tr}
			\left[(\log\mathbf{D}_N)^2\right].
		\end{equation}
		Combining this result with \eqref{eq:log_frobenius_convergence}
		yields
		\begin{equation}
			\begin{aligned}
				\lim_{N\to\infty}
				\frac{1}{N}
				d_{\mathrm{AIRM}}^2(\mathbf{A}_N,\mathbf{B}_N)
				&=
				\lim_{N\to\infty}
				\frac{1}{N}
				\|\log\mathbf{H}_N\|_{\mathrm F}^2\\
				&=
				\frac{1}{2\pi}
				\int_{-\pi}^{\pi}
				\left|
				\log S_2(\omega)-\log S_1(\omega)
				\right|^2\,d\omega\\
				&=
				d_{\mathrm{FR}}^2(S_1,S_2).
			\end{aligned}
		\end{equation}
		Taking the square root completes the proof.
		
	\end{proof}
	\section{Extension to Elliptical Observation Models}
	\label{sec:elliptical_extension}
	
	The goal of this section is to explore how far this intrinsic representation can extend beyond the Gaussian WSS setting.
	
	To this end, we consider finite-dimensional elliptical observation models whose covariance matrices are induced by WSS processes. This framework accommodates a range of non-Gaussian distributions, including K-distributed and Student's t-distributed observations, while retaining a covariance-based description of second-order dependence.

	\subsection{Elliptical Observation Models}
	\label{subsec:elliptical_model}
	
	The Gaussian construction in the previous section provides an intrinsic geometric description of stationary covariance structures. We now investigate how this representation can be extended beyond Gaussian observations, particularly to finite-dimensional models exhibiting heavier-tailed distributions or random variations in signal power.
	
	Consider two independent realizations of a common-scale Gaussian mixture process, and let \(\tau_i\) denote the random scale associated with realization \(i\). Within each realization, the same random scale is shared by the entire process, yielding the conditional covariance matrices
	\[
	R_{N,i}^{\mathrm{cond}} = \tau_i R_N,
	\qquad i=1,2,
	\]
	even though both realizations are generated from the same underlying covariance structure \(R_N\). Such variations introduce differences in overall power without necessarily changing the underlying covariance shape. To focus on the structural differences of interest, it is therefore natural to disentangle the global scale from the covariance shape.

	This distinction is useful in signal-processing applications where overall power may vary across realizations or observation conditions, while the relative covariance structure remains the object of interest. A geometric representation that separates global scale from covariance shape provides a focused description of these structural variations.

	These considerations motivate us to investigate the intrinsic geometry of covariance shape within the framework of elliptically distributed finite-dimensional observation models \cite{cambanis1981theory,fang2018symmetric,ollila2012complex}.
	
	Let $S$ be a strictly positive spectral density satisfying the
	assumptions in Section 3, and let
	\[
	R_N(S) = T_N(S)
	\]
	be the covariance matrix induced by an observation window of
	length $N$. We consider a zero-mean proper complex elliptical
	observation vector $y_N \in \mathbb{C}^N$ with density
	\[
	p_{N,S,\theta}(y_N)
	=
	\frac{1}{\det R_N(S)}
	h_{N,\theta}\left(
	y_N^H R_N(S)^{-1} y_N
	\right),
	\]
	where $h_{N,\theta}$ is a radial generator parameterized by
	$\theta$. It is normalized so that the density integrates to one. We assume the regularity conditions required for the
	corresponding Fisher information to be well-defined,
	which will be specified further in Section \ref{subsec:radial_condition}. 
	We additionally require
	\[
	\mathbb{E}[q_N] = N,
	\qquad
	q_N = y_N^H R_N(S)^{-1} y_N,
	\]
	so that $R_N(S)$ is the covariance matrix of the observation.
	
	The finite-dimensional observation models are constructed from
	the covariance structure of an underlying infinite-length
	stationary process through finite observation windows.
	If they are to represent the finite-dimensional marginals
	of a single infinite-length process, they must satisfy
	projective consistency. A projectively consistent example
	will be provided in Section~\ref{subsec:radial_condition}.
	
	For convenience, we decompose the covariance matrix as
	\[
	R_N = \beta_N C_N,
	\qquad
	\beta_N > 0,
	\qquad
	\det C_N = 1,
	\]
	where
	\[
	\beta_N = (\det R_N)^{1/N}
	\]
	represents the global scale and $C_N$ the normalized shape.
	Define
	\[
	q = y_N^H R_N^{-1}y_N,
	\qquad
	u = \frac{R_N^{-1/2}y_N}{\sqrt{q}},
	\qquad
	\phi_{N,\theta}(q)
	=
	-q\frac{\partial}{\partial q}\log h_{N,\theta}(q).
	\]
	By the spherical symmetry of the whitened elliptical distribution \cite{kelker1970distribution,cambanis1981theory,ollila2012complex},
	$u$ is uniformly distributed on the complex unit sphere and
	independent of the squared radius $q$.
	
\begin{figure*}[h]
		\centering
		\includegraphics[width=0.85\textwidth]
		{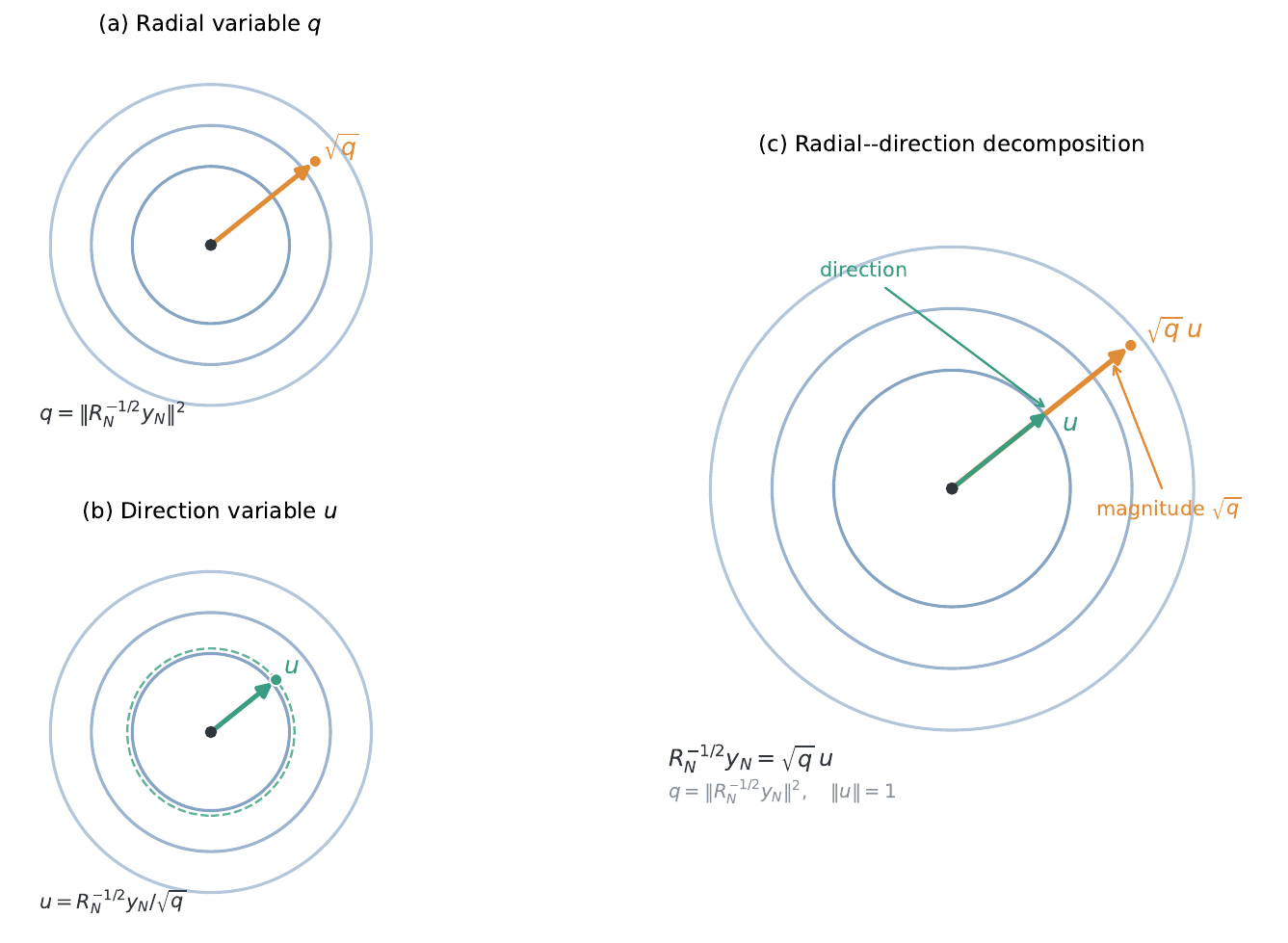}
		\caption{Radial--direction decomposition of a whitened
			elliptical observation. (a) The squared Mahalanobis radius
			$q$ determines the radial magnitude. (b) The normalized
			direction $u$ lies on the complex unit sphere, schematically
			represented in two dimensions. (c) Their correspondence
			through $R_N^{-1/2}y_N=\sqrt{q}\,u$.}
		\label{fig:radial_direction}
	\end{figure*}

	Figure~\ref{fig:radial_direction} illustrates the radial--direction
	decomposition of a whitened elliptical observation.
	
	\subsection{Fisher Information and Orthogonal Decomposition}
	\label{subsec:elliptical_fisher}
	
	Let $\dot C_N$ be a tangent direction at $C_N$, and define its
	whitened representation as
	\begin{equation}
		A=C_N^{-1/2}\dot C_N C_N^{-1/2},
		\qquad \operatorname{tr}(A)=0,
	\end{equation}
	where the trace constraint follows from $\det C_N=1$.
	
	With the global scale $\beta_N$ held fixed, the corresponding
	covariance perturbation is $\dot R_N=\beta_N\dot C_N$.
	The log-determinant term has zero derivative:
	\begin{equation}
		\operatorname{tr}(R_N^{-1}\dot R_N)
		=\operatorname{tr}(C_N^{-1}\dot C_N)
		=\operatorname{tr}(A)=0.
	\end{equation}
	Meanwhile, differentiating the quadratic form gives
	\begin{equation}
		\dot q
		=-y_N^H R_N^{-1}\dot R_N R_N^{-1}y_N
		=-q\,u^H A u.
	\end{equation}
	Consequently, the shape score is
	\begin{equation}
		s_A(y_N)
		=\phi_{N,\theta}(q)\,u^H A u.
		\label{eq:elliptical_shape_score}
	\end{equation}
	
	For Hermitian trace-free matrices $A$ and $B$, the complex-sphere
	moment identity gives \cite{collins2006integration}
	\begin{equation}
		\mathbb E\big[(u^H A u)(u^H B u)\big]
		=\frac{\operatorname{tr}(AB)}{N(N+1)}.
		\label{eq:complex_sphere_moment}
	\end{equation}
	Using the independence of $u$ and $q$, we obtain
	\begin{align}
		g_{\mathrm{shape},\mathrm{shape};N,\theta}(A,B)
		&=\mathbb E[s_A(y_N)s_B(y_N)]\\
		&=\mathbb E[\phi_{N,\theta}(q)^2]\,
		\mathbb E[(u^H A u)(u^H B u)]\\
		&=c_{N,\theta}\operatorname{tr}(AB),
		\label{eq:elliptical_shape_metric}
	\end{align}
	where
	\begin{equation}
		c_{N,\theta}
		=\frac{\mathbb E[\phi_{N,\theta}(q)^2]}{N(N+1)}.
	\end{equation}
	
	We next examine the cross-information between shape and nuisance
	directions. Since $u$ and $q$ are independent and
	\[
	\mathbb E[u^H A u]
	=\frac{\operatorname{tr}(A)}{N}=0,
	\]
	the shape score is orthogonal to both the global scale score and
	the radial-parameter scores. Consequently,
	\begin{equation}
		g_{\mathrm{shape},\mathrm{scale};N,\theta}=0,
		\qquad
		g_{\mathrm{shape},\mathrm{radial};N,\theta}=0.
		\label{eq:shape_radial_orthogonality}
	\end{equation}
	
\begin{figure}[h]
		\centering
		\includegraphics[width=0.7\linewidth]
		{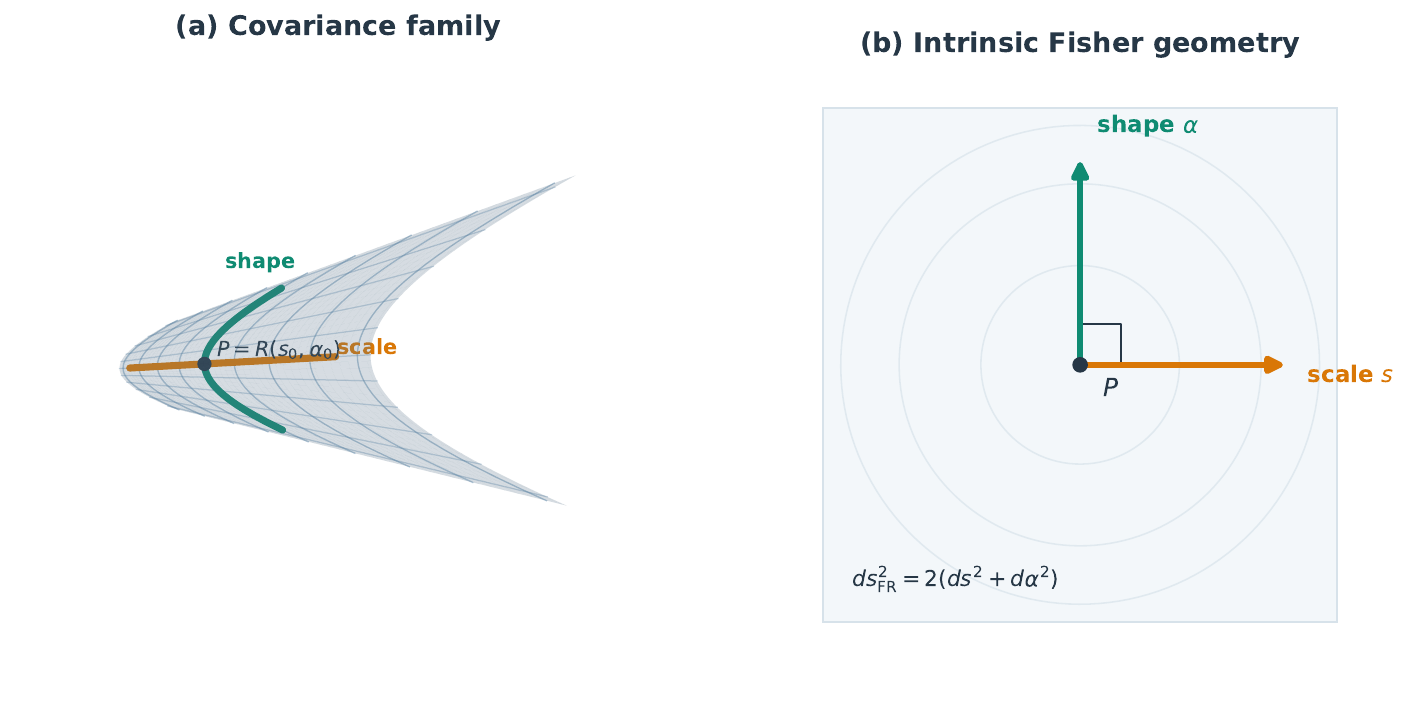}
		\caption{An example of scale--shape decomposition in a
			two-dimensional covariance family. (a) Scale and shape
			variations in matrix-entry coordinates. (b) Their
			orthogonality in the intrinsic Fisher geometry.}
		\label{fig:scale_shape}
	\end{figure}

	Figure~\ref{fig:scale_shape} illustrates an example of
	the orthogonal scale--shape decomposition.
	
	For a fixed radial family and fixed $\theta$, the coefficient
	$c_{N,\theta}$ is constant over the determinant-one shape manifold.
	The induced finite-dimensional shape distance is therefore
	\begin{equation}
		d_{\mathrm{FR,shape},N,\theta}(C_1,C_2)
		=\sqrt{c_{N,\theta}}\,
		\left\|\log\left(C_1^{-1/2}C_2C_1^{-1/2}\right)\right\|_F.
		\label{eq:finite_shape_distance}
	\end{equation}
	If $\theta$ varies along a path, the shape contribution to its
	length is instead given by
	\begin{equation}
		L_{\mathrm{shape}}
		=\int_0^1
		\sqrt{c_{N,\theta(t)}
			\operatorname{tr}\big(A(t)^2\big)}\,dt,
	\end{equation}
	where
	\[
	A(t)=C_N(t)^{-1/2}\dot C_N(t)C_N(t)^{-1/2}.
	\]
	
\subsection{Gaussian Shape Recovery Condition}
	\label{subsec:radial_condition}
	
	Recall that the squared radius $q=y_N^H R_N^{-1}y_N$ has density
	\begin{equation}
		f_{Q,N,\theta}(q)
		=K_{N,\theta}q^{N-1}h_{N,\theta}(q),
		\qquad q>0,
	\end{equation}
	where $K_{N,\theta}$ is a normalization constant. The factor
	$q^{N-1}$ arises from the radial Jacobian in the $N$-dimensional
	complex space, while $h_{N,\theta}(q)$ describes the radial
	distribution of the elliptical model.
	
	To characterize the radial contribution to the shape Fisher
	information, it is convenient to introduce the log-radius
	$T=\log q$. By the change of variables, its density is
	\begin{equation}
		f_{T,N,\theta}(t)
		=e^t f_{Q,N,\theta}(e^t)
		\propto e^{Nt}h_{N,\theta}(e^t).
	\end{equation}
	Assume that $f_{T,N,\theta}$ is differentiable, vanishes at
	both boundaries of its support, and satisfies the integrability
	conditions required for integration by parts. Then the expected
	log-density derivative is zero, i.e.,
	\begin{equation}
		\mathbb E\left[
		\frac{d}{dt}\log f_{T,N,\theta}(T)
		\right]=0.
	\end{equation}
	Consequently,
	\begin{equation}
		\frac{d}{dt}\log f_{T,N,\theta}(t)
		=N-\phi_{N,\theta}(e^t),
		\qquad
		\mathbb E[\phi_{N,\theta}(q)]=N.
		\label{eq:radial_score_identity}
	\end{equation}
	
	Define the location Fisher information of $T=\log q$ by
	\begin{equation}
		J_{N,\theta}
		=\mathbb E\left[
		\left(
		\frac{d}{dt}\log f_{T,N,\theta}(T)
		\right)^2
		\right].
		\label{eq:log_radial_fisher}
	\end{equation}
	It follows from \eqref{eq:radial_score_identity} that
	\begin{equation}
		J_{N,\theta}
		=\operatorname{Var}(\phi_{N,\theta}(q)),
	\end{equation}
	and hence
	\begin{equation}
		c_{N,\theta}
		=\frac{N}{N+1}
		+\frac{J_{N,\theta}}{N(N+1)}.
		\label{eq:shape_coefficient_identity}
	\end{equation}
	
	We therefore obtain the following criterion.
	
	\begin{proposition}[Radial condition for Gaussian shape-geometry recovery]
		\label{prop:Radial_cond}
		Suppose the regularity and boundary assumptions above hold. Then
		\[
		c_{N,\theta}\to 1
		\quad\Longleftrightarrow\quad
		J_{N,\theta}=o(N^2).
		\]
		
		If convergence is required uniformly over a collection of radial
		parameters, the corresponding necessary and sufficient uniform
		condition is
		\begin{equation}
			\sup_{\theta\in\Theta}J_{N,\theta}=o(N^2).
		\end{equation}
	\end{proposition}
	The condition concerns the radial score, not the tail label of
	the distribution. In particular, heavy tails alone do not prevent
	the recovery of the Gaussian shape geometry.
	
	A broad sufficient class is provided by common-scale Gaussian
	mixtures \cite{andrews1974scale,shah2014student}, which also gives a projectively consistent
	infinite-dimensional process model. Consider
	\begin{equation}
		z[n]=\sqrt{\tau}\,x[n],
	\end{equation}
	where $\tau>0$ is a random variable independent of the entire
	zero-mean proper complex Gaussian WSS process $x[n]$, with
	$\mathbb E[\tau]=1$. The same random scale $\tau$ is shared
	by the entire process, rather than independently drawn for each
	observation window. For an observation window of length $N$,
	the marginal density of $z_N$ is
	\begin{equation}
		p(z_N)
		=\frac{1}{\det R_N}h_N(q),
		\qquad q=z_N^H R_N^{-1}z_N,
	\end{equation}
	where
	\begin{equation}
		h_N(q)
		=\frac{1}{\pi^N}
		\mathbb E_\tau\left[
		\tau^{-N}\exp\left(-\frac{q}{\tau}\right)
		\right].
		\label{eq:mixture_radial_generator}
	\end{equation}
	Thus, the radial generators $h_N$ are induced by a single,
	dimension-independent mixing distribution. Since all finite-dimensional
	distributions are marginals of the same process, they automatically
	satisfy projective consistency. The corresponding radial score is
	\begin{equation}
		\phi_N(q)
		=q\frac{
			\mathbb E_\tau[
			\tau^{-N-1}\exp(-q/\tau)
			]
		}{
			\mathbb E_\tau[
			\tau^{-N}\exp(-q/\tau)
			]
		}.
		\label{eq:mixture_radial_score}
	\end{equation}
	
	Now let $V=x_N^H R_N^{-1}x_N$. Since $x_N$ is proper complex
	Gaussian, $V\sim\operatorname{Gamma}(N,1)$ and $q=\tau V$.
	The posterior density of $\tau$ given $q$ is
	\begin{equation}
		p(\tau\mid q)
		=
		\frac{
			p_\tau(\tau)\tau^{-N}\exp(-q/\tau)
		}{
			\mathbb E_\tau[
			\tau^{-N}\exp(-q/\tau)
			]
		}.
	\end{equation}
	Consequently,
	\begin{equation}
		\begin{aligned}
			\mathbb E[V\mid q]
			&=q\mathbb E\left[
			\frac{1}{\tau}\middle|q
			\right]\\
			&=q\frac{
				\mathbb E_\tau[
				\tau^{-N-1}\exp(-q/\tau)
				]
			}{
				\mathbb E_\tau[
				\tau^{-N}\exp(-q/\tau)
				]
			}\\
			&=\phi_N(q).
		\end{aligned}
		\label{eq:mixture_conditional_score}
	\end{equation}
	
	The conditional-variance identity therefore gives
	\begin{equation}
		\begin{aligned}
			\operatorname{Var}(\phi_N(q))
			&=\operatorname{Var}(\mathbb E[V\mid q])\\
			&\leq\operatorname{Var}(V)\\
			&=N.
		\end{aligned}
	\end{equation}
	Combining this bound with \eqref{eq:shape_coefficient_identity},
	we obtain
	\begin{equation}
		\frac{N}{N+1}
		\leq c_N\leq 1,
		\label{eq:mixture_coefficient_bound}
	\end{equation}
	and hence
	\begin{equation}
		0\leq 1-c_N\leq\frac{1}{N+1}.
	\end{equation}
	This bound guarantees that $c_N\to1$ as $N\to\infty$.
	In particular, this projectively consistent non-Gaussian
	process model satisfies the radial Fisher-information condition
	required for asymptotic Gaussian shape-geometry recovery.
	
	\subsection{Asymptotic Spectral-Shape Geometry}
	\label{subsec:elliptical_spectral_rate}
	
	We now specialize to Toeplitz covariance matrices
	$R_N(S)=T_N(S)$. Following the scale-shape decomposition
	introduced in Section \ref{subsec:elliptical_model}, define
	\begin{equation}
		d_N(S)=\det(R_N(S))^{1/N},
		\qquad
		C_N(S)=\frac{R_N(S)}{d_N(S)}.
		\label{eq:normalized_toeplitz_covariance}
	\end{equation}
	Thus, $d_N(S)$ represents the global scale and
	$\det C_N(S)=1$.
	
	For a spectral perturbation $H$, define
	\begin{equation}
		B_N(H)
		=R_N(S)^{-1/2}T_N(H)R_N(S)^{-1/2}.
	\end{equation}
	Differentiating the determinant normalization in
	\eqref{eq:normalized_toeplitz_covariance} gives the
	corresponding trace-free shape tangent
	\begin{equation}
		A_N(H)
		=B_N(H)-\frac{1}{N}\operatorname{tr}(B_N(H))I_N,
		\qquad
		\operatorname{tr}(A_N(H))=0.
		\label{eq:toeplitz_shape_tangent}
	\end{equation}
	Substituting these tangent vectors into
	\eqref{eq:elliptical_shape_metric} yields
	\begin{align}
		&g_{\mathrm{shape},\mathrm{shape};N,\theta}(H,K)
		=c_{N,\theta}\operatorname{tr}
		\big(A_N(H)A_N(K)\big) \nonumber\\
		&=c_{N,\theta}\Bigg[
		\operatorname{tr}\big(B_N(H)B_N(K)\big)
		-\frac{1}{N}
		\operatorname{tr}(B_N(H))
		\operatorname{tr}(B_N(K))
		\Bigg].
		\label{eq:toeplitz_shape_fisher}
	\end{align}
	
	Let $H=S\dot u$ and $K=S\dot v$, where $\dot u$ and
	$\dot v$ are continuous real-valued perturbations of the
	log-spectrum. Under the same Toeplitz regularity assumptions
	used in Section \ref{sec:gaussian_model}, the normalized trace limits are
	\begin{align}
		\lim_{N\to\infty}\frac{1}{N}
		\operatorname{tr}\big(B_N(H)B_N(K)\big)
		&=\frac{1}{2\pi}\int_{-\pi}^{\pi}
		\dot u(\omega)\dot v(\omega)\,d\omega,\\
		\lim_{N\to\infty}\frac{1}{N}
		\operatorname{tr}(B_N(H))
		&=\frac{1}{2\pi}\int_{-\pi}^{\pi}
		\dot u(\omega)\,d\omega.
	\end{align}
	For notational convenience, define the frequency average
	\begin{equation}
		\overline{f}
		=\frac{1}{2\pi}\int_{-\pi}^{\pi}f(\omega)\,d\omega.
	\end{equation}
	
	The preceding calculations yield the following result.
	
	\begin{theorem}[Asymptotic spectral-shape Fisher information rate]
		\label{thm:spectral_shape_rate}
		Let $S$ be continuous and uniformly bounded away from zero,
		and let $H,K$ be continuous real-valued perturbations.
		Suppose that the finite-window elliptical models satisfy
		the radial regularity assumptions of Proposition \ref{prop:Radial_cond} and
		\begin{equation}
			J_{N,\theta}=o(N^2).
		\end{equation}
		Then, for $H=S\dot u$ and $K=S\dot v$,
		\begin{equation}
			\lim_{N\to\infty}\frac{1}{N}
			g_{\mathrm{shape},\mathrm{shape};N,\theta}(H,K)
			=
			\frac{1}{2\pi}\int_{-\pi}^{\pi}
			(\dot u-\overline{\dot u})
			(\dot v-\overline{\dot v})\,d\omega.
			\label{eq:elliptical_spectral_shape_metric}
		\end{equation}
		If
		\begin{equation}
			\sup_{\theta\in\Theta}J_{N,\theta}=o(N^2),
		\end{equation}
		the convergence is uniform over $\theta\in\Theta$.
	\end{theorem}
	The result follows directly by substituting the trace limits and applying
	Proposition~\ref{prop:Radial_cond}; the uniform statement follows in the same way.
	
	Equation~\eqref{eq:elliptical_spectral_shape_metric}
	defines the intrinsic Fisher metric on log-spectral shapes
	modulo global positive scaling.
	
	\subsection{Geodesic Distance and Relation to the Gaussian Construction}
	\label{subsec:elliptical_shape_distance}

	Define the centered log-spectrum coordinate
	\begin{equation}
		u_c(\omega)
		=\log S(\omega)
		-\frac{1}{2\pi}\int_{-\pi}^{\pi}
		\log S(\xi)\,d\xi.
	\end{equation}
	In this coordinate, the limiting metric takes the form
	\begin{equation}
		g(\dot u_c,\dot v_c)
		=\frac{1}{2\pi}\int_{-\pi}^{\pi}
		\dot u_c(\omega)\dot v_c(\omega)\,d\omega,
	\end{equation}
	which is the ordinary $L^2$ inner product on the space
	of zero-mean log-spectra.
	
	Since this metric is independent of the base point $u_c$,
	the induced geometry is flat. Consequently, the minimizing
	geodesic between two spectral shapes is given by linear
	interpolation in centered log-spectrum coordinates.
	The resulting geodesic distance is therefore
	\begin{equation}
		d_{\mathrm{FR,shape}}^2(S_1,S_2)
		=\frac1{2\pi}\int_{-\pi}^{\pi}
		\left[
		\log\frac{S_2(\omega)}{S_1(\omega)}
		-\overline{\log(S_2/S_1)}
		\right]^2d\omega,
		\label{eq:elliptical_spectral_shape_distance}
	\end{equation}
	where the overline denotes the frequency average.
	This distance is invariant under independent positive
	rescalings of either spectrum:
	\begin{equation}
		d_{\mathrm{FR,shape}}(aS_1,bS_2)
		=d_{\mathrm{FR,shape}}(S_1,S_2),
		\qquad a,b>0.
	\end{equation}
	Equivalently, it is the quotient distance induced by the
	Gaussian spectral Fisher--Rao metric after identifying
	spectra that differ only by a positive constant multiplier:
	\begin{equation}
		d_{\mathrm{FR,shape}}([S_1],[S_2])
		=\inf_{a>0}d_{\mathrm{FR}}(S_1,aS_2).
	\end{equation}
	
	The determinant normalization removes the global scale component
	of the finite-window covariance geometry. To see its effect on the
	AIRM distance, let
	\[
	H_N=R_N(S_1)^{-1/2}R_N(S_2)R_N(S_1)^{-1/2},
	\qquad
	\alpha_N=\log\frac{d_N(S_1)}{d_N(S_2)}.
	\]
	The relative matrix formed from the normalized covariances is
	then $(d_N(S_1)/d_N(S_2))H_N$, so that
	\begin{align}
		\frac1N d_{\mathrm{AIRM}}^2(C_N(S_1),C_N(S_2))
		&=\frac1N\left\|\log H_N+\alpha_N I_N\right\|_{\mathrm F}^2\nonumber\\
		&=\frac1N\operatorname{tr}\big[(\log H_N)^2\big]
		-\left(\frac1N\operatorname{tr}(\log H_N)\right)^2,
		\label{eq:normalized_airm_identity}
	\end{align}
	where the second equality follows from
	$\alpha_N=-N^{-1}\operatorname{tr}(\log H_N)$.
	Theorem~\ref{thm:airm_consistency} and the Toeplitz approximation established
	in its proof yield the limit of the first term, while
	the same approximation together with the first Szegő
	limit theorem gives
	\[
	\frac1N\operatorname{tr}(\log H_N)
	\longrightarrow
	\overline{\log(S_2/S_1)}.
	\]
	Therefore,
	\begin{equation}
		\lim_{N\to\infty}\frac1N
		d_{\mathrm{AIRM}}^2(C_N(S_1),C_N(S_2))
		=d_{\mathrm{FR,shape}}^2(S_1,S_2).
		\label{eq:normalized_airm_limit}
	\end{equation}
	Thus, the spectral-shape distance coincides with the
	asymptotic normalized AIRM distance between
	determinant-normalized covariance representations.
	
	\subsection{Discussion and Scope}
	\label{subsec:elliptical_scope}
	
	The results extend the spectral-domain information geometry
	developed for Gaussian WSS processes to elliptically
	distributed observations. Under the radial
	Fisher-information condition $J_{N,\theta}=o(N^2)$,
	the limiting intrinsic distance between two spectral
	shapes is given by
	\[
	d_{\mathrm{FR,shape}}(S_1,S_2)
	=
	\left\{
	\frac{1}{2\pi}\int_{-\pi}^{\pi}
	\left[
	\log\frac{S_2(\omega)}{S_1(\omega)}
	-
	\overline{\log\frac{S_2}{S_1}}
	\right]^2
	d\omega
	\right\}^{1/2}.
	\]
	This coincides with the Gaussian spectral
	Fisher--Rao distance after quotienting out the global
	power scale. This shows that the intrinsic spectral-shape geometry
	developed in the Gaussian setting extends to a broad
	class of elliptically distributed observation models.
	
	\section{Numerical Experiments}
	\label{sec:numerical}
	
	This section numerically examines the main asymptotic results developed
	in Sections~3 and~4. The experiments address three aspects:
	the Gaussian Fisher information rate, the consistency between
	finite-window AIRM distances and their spectral Fisher--Rao limits,
	the recovery of the Gaussian shape coefficient for projectively consistent non-Gaussian examples and its behavior under critical radial Fisher-information growth.
	 These experiments
	illustrate the theoretical results and are not intended as
	substitutes for their proofs.
	
	\subsection{Gaussian Fisher Information Rate}
	\label{subsec:exp_fisher_rate}
	
\begin{figure}[h]
		\centering
		\includegraphics[width=0.8\linewidth]
		{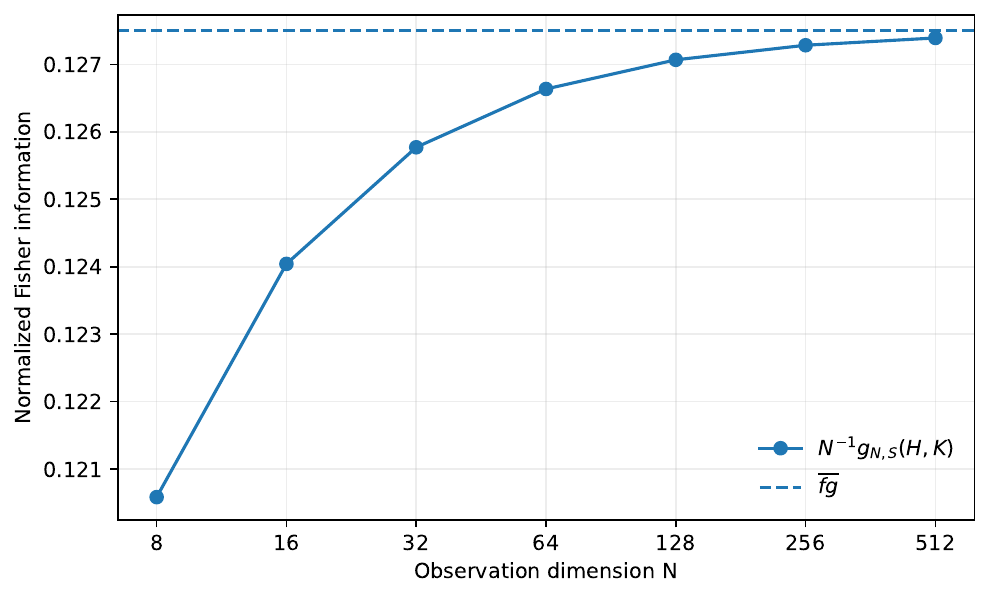}
		\caption{Convergence of the normalized finite-dimensional
			Gaussian Fisher information to its spectral limit.}
		\label{fig:exp_fisher_rate}
	\end{figure}

	We first examine the convergence of the finite-dimensional Gaussian
	Fisher information rate. Consider the strictly positive spectral
	density
	\[
	S(\omega)=\exp\!\big(0.55\cos\omega+0.22\cos 2\omega\big),
	\]
	and the relative perturbations
	\[
	f(\omega)=0.70\cos\omega-0.25\sin 2\omega,
	\qquad
	g(\omega)=0.40\cos\omega+0.30\cos 2\omega+0.10\sin 2\omega.
	\]
	We set $H=Sf$ and $K=Sg$. For each observation dimension $N$, the
	finite-window Fisher information is computed from
	\[
	g_{N,S}(H,K)
	=\operatorname{tr}\!\left[
	R_N(S)^{-1}T_N(H)R_N(S)^{-1}T_N(K)
	\right].
	\]
	The limiting value predicted by the Fisher information rate is
	\[
	\overline{fg}
	=\frac{1}{2\pi}\int_{-\pi}^{\pi}
	f(\omega)g(\omega)\,d\omega
	=0.1275.
	\]
	Figure~\ref{fig:exp_fisher_rate} shows that
	$N^{-1}g_{N,S}(H,K)$ approaches this value as the observation
	dimension increases. The remaining discrepancy decreases
	monotonically over the tested dimensions.
	
	\subsection{Finite-Window AIRM and Spectral Fisher--Rao Distances}
	\label{subsec:exp_airm}
	
\begin{figure*}[h]
		\centering
		\includegraphics[width=0.96\textwidth]
		{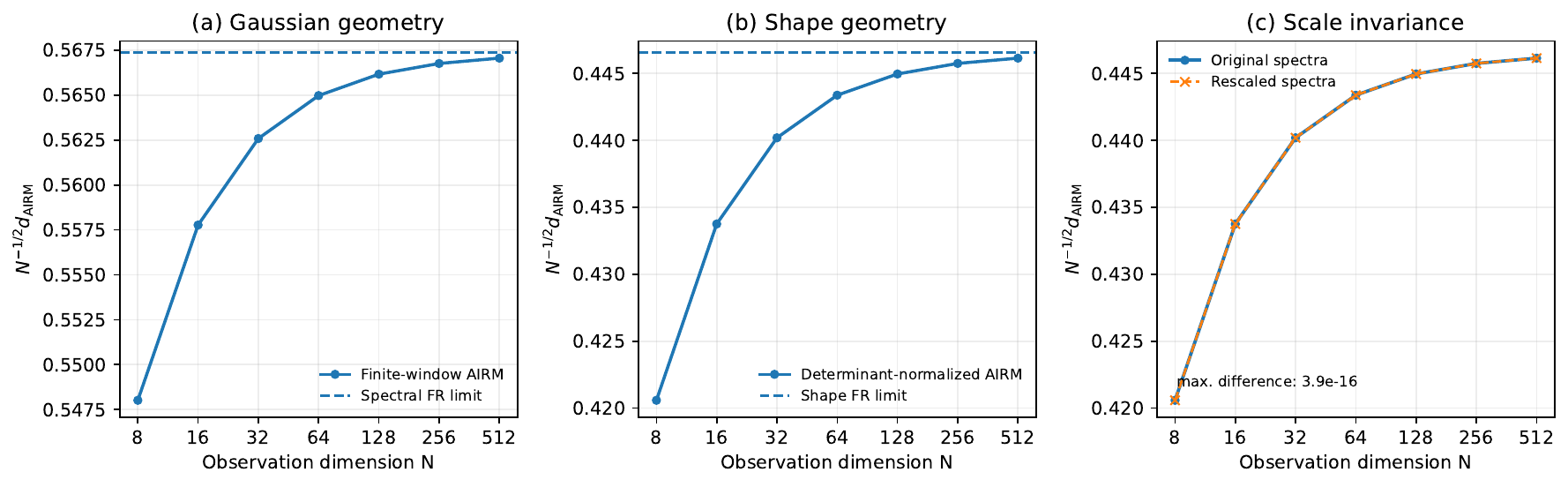}
		\caption{Finite-window AIRM convergence. (a) The normalized
			AIRM approaches the Gaussian spectral Fisher--Rao distance.
			(b) The determinant-normalized AIRM approaches the
			spectral-shape distance. (c) Independent positive rescalings
			leave the determinant-normalized AIRM unchanged up to
			numerical precision.}
		\label{fig:exp_airm}
	\end{figure*}

	We next compare the finite-window AIRM with the Gaussian spectral
	Fisher--Rao distance and its scale-invariant shape counterpart.
	Let
	\[
	S_1(\omega)=e^{u_1(\omega)},\qquad
	S_2(\omega)=e^{u_2(\omega)},
	\]
	where
	\[
	u_1(\omega)=0.42\cos\omega+0.16\cos 2\omega,
	\]
	\[
	u_2(\omega)=0.35-0.18\cos\omega
	+0.34\cos 2\omega+0.08\sin\omega.
	\]
	The Gaussian spectral distance is
	\[
	d_{\mathrm{FR}}(S_1,S_2)
	=\sqrt{\overline{(u_2-u_1)^2}},
	\]
	whereas the spectral-shape distance is
	\[
	d_{\mathrm{FR,shape}}(S_1,S_2)
	=\sqrt{\overline{\left[
			(u_2-u_1)-\overline{(u_2-u_1)}
			\right]^2}}.
	\]
	For each $N$, we compute the AIRM distances between
	$R_N(S_1)$ and $R_N(S_2)$, and between their determinant-normalized
	versions
	\[
	C_N(S_i)=\frac{R_N(S_i)}{\det(R_N(S_i))^{1/N}}.
	\]
	Figure~\ref{fig:exp_airm} compares the normalized finite-window
	distances with their spectral limits. Both sequences approach
	the corresponding limits as $N$ increases. The nonzero mean
	$\overline{u_2-u_1}=0.35$ makes the distinction between the
	full Gaussian spectral distance and the shape distance explicit:
	the latter removes the global log-scale difference.
	
	We additionally replace the two spectra by $3.7S_1$ and
	$0.23S_2$ before forming the determinant-normalized covariance
	matrices. The two resulting shape-AIRM sequences overlap to numerical
	precision; their maximum discrepancy over the tested dimensions
	is below $2\times 10^{-16}$. This confirms the invariance to
	independent positive rescalings.
	
	 \subsection{Convergence of the Elliptical Shape Coefficient}
	\label{subsec:exp_radial}
	
	Finally, we investigate the radial coefficient
	\[
	c_N=\frac{\mathbb{E}[\phi_N(q)^2]}{N(N+1)}
	=\frac{N}{N+1}+\frac{J_N}{N(N+1)}.
	\]
	We first consider the common-scale Gaussian mixture models discussed in Section~4.3, including the Gaussian baseline and two non-Gaussian mixing distributions: a Gamma texture
	\(\tau\sim\operatorname{Gamma}(2,\text{rate}=2)\) and an inverse-Gamma texture
	\(\tau\sim\operatorname{InvGamma}(3,\text{scale}=2)\).
	Both mixing distributions have unit mean. In each non-Gaussian model, a single random scale \(\tau\) is shared by the entire underlying Gaussian process, ensuring that the resulting finite-dimensional distributions are projectively consistent.
	
	For each model and dimension, the coefficient is estimated by Monte Carlo sampling of \(q=\tau V\), where \(V\sim\operatorname{Gamma}(N,1)\), and evaluating the corresponding radial score. We use \(250{,}000\) Monte Carlo samples for each model and dimension. The Gaussian value is exactly \(c_N=1\). For both non-Gaussian mixtures, the estimates approach one as \(N\) increases. The empirical deficits also remain below the theoretical upper bound \(1/(N+1)\) over all tested dimensions. At \(N=128\), for example, the estimated coefficients are \(0.992367\) and \(0.992424\) for the Gamma and inverse-Gamma mixtures, respectively. The normalized empirical means of the radial scores are close to one, consistent with \(\mathbb{E}[\phi_N(q)]=N\).
	
	To further examine the role of the radial Fisher-information growth condition, we introduce a dimension-dependent radial model as a contrasting example. Specifically, let the log-radius satisfy
	\[
	T_N=\log q_N\sim\mathcal{N}(\mu_N,\sigma_N^2),
	\qquad
	\sigma_N^2=\frac{1}{\kappa N^2},
	\]
	where
	\[
	\mu_N=\log N-\frac{1}{2\kappa N^2},
	\qquad \kappa>0.
	\]
	This choice ensures that \(\mathbb{E}[q_N]=N\), while the location Fisher information of \(T_N\) is
	\[
	J_N=\frac{1}{\sigma_N^2}=\kappa N^2.
	\]
	Consequently, the shape coefficient is given by
	\[
	c_N=\frac{N(1+\kappa)}{N+1}
	\longrightarrow 1+\kappa,
	\]
	rather than converging to its Gaussian value. We use \(\kappa=0.5\) in the numerical experiment. Unlike the common-scale Gaussian mixtures, this dimension-dependent family is introduced solely as a finite-dimensional counterexample and is not required to be projectively consistent.
	
\begin{figure}[h]
		\centering
		\includegraphics[width=0.7\linewidth]
		{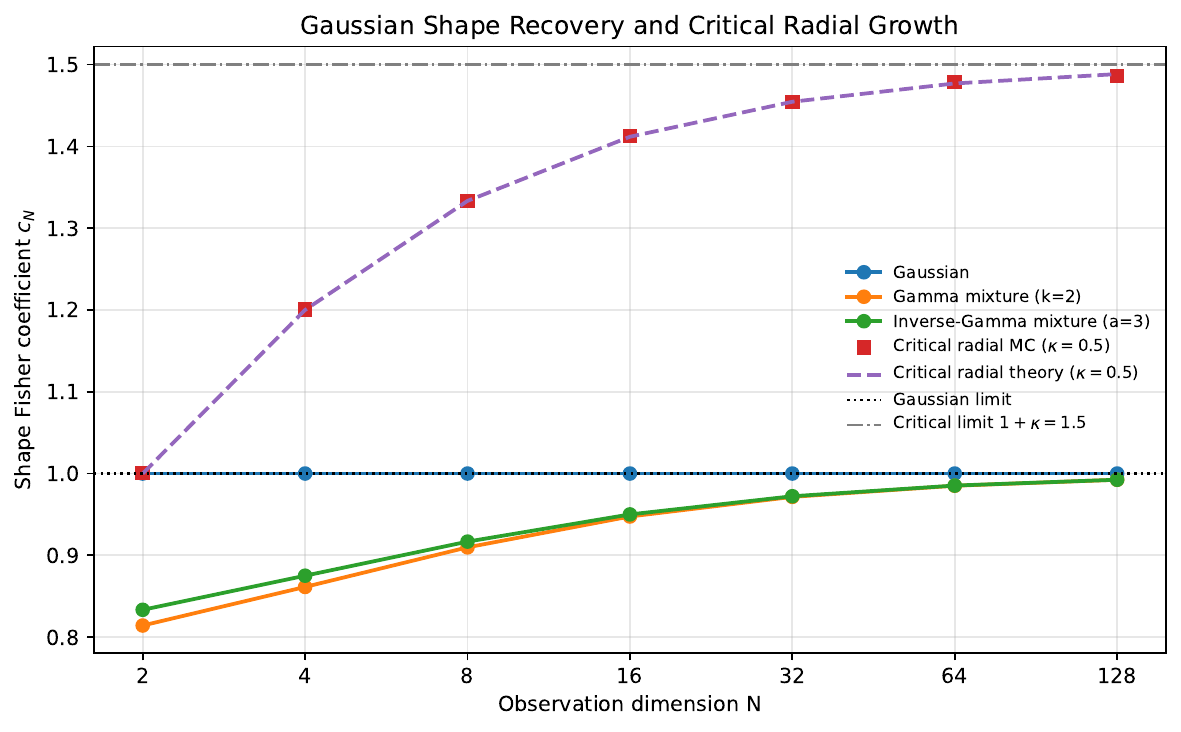}
		\caption{Radial Fisher-information growth and the elliptical shape coefficient. The common-scale Gaussian mixtures exhibit convergence of \(c_N\) to the Gaussian value one. In contrast, the dimension-dependent critical radial model with \(J_N=0.5N^2\) converges to \(1.5\), illustrating the effect of quadratic radial Fisher-information growth.}
		\label{fig:exp_radial}
	\end{figure}

	Figure~\ref{fig:exp_radial} compares the Gaussian recovery of the common-scale mixtures with the non-Gaussian limiting behavior of the critical radial model. These results illustrate the role of the condition \(J_N=o(N^2)\): subquadratic radial Fisher-information growth leads to recovery of the Gaussian shape coefficient, whereas quadratic growth can produce a non-Gaussian limiting coefficient.
	
	\section{Conclusion}
	\label{sec:conclusion}
	
	This work develops an information-geometric framework for stationary
	covariance structures through finite-dimensional observation models.
	For proper complex Gaussian WSS processes, the normalized
	finite-window Fisher information converges to an $L^2$ metric on
	log-spectral perturbations, and the corresponding spectral
	Fisher--Rao distance is recovered as the asymptotic AIRM distance
	between Toeplitz covariance matrices. The construction is then
	extended to elliptically distributed observations by separating
	global scale from covariance shape. Under a subquadratic radial
	Fisher-information condition, the normalized shape Fisher
	information converges to the centered $L^2$ metric on log-spectral
	shapes, yielding a scale-invariant spectral-shape distance.
	Common-scale Gaussian mixtures provide a projectively consistent
	non-Gaussian class satisfying this condition. A dimension-dependent
	critical radial model further illustrates that quadratic growth of
	the radial Fisher information can produce a non-Gaussian limiting
	shape coefficient. The numerical experiments illustrate the Fisher
	information rate, the two AIRM limits, and both behaviors of the
	elliptical shape coefficient. Together, these results provide a unified geometric
	description of stationary covariance and spectral-shape structures
	for Gaussian and a broad class of elliptical observation models.
	
	\bibliographystyle{unsrtnat}
	\bibliography{references_arxiv}
	
\end{document}